\documentclass [12pt,a4paper]{article}
\usepackage[T1]{fontenc}
\usepackage[ansinew]{inputenc}
\usepackage{amssymb}
\usepackage{amsmath}
\usepackage{graphicx}
\usepackage{amsthm}
\usepackage{cite}

\newtheorem{theorem}{\textbf{Theorem}}[section]
\newtheorem{lemma}{\textbf{Lemma}}[section]
\newtheorem{proposition}{\textbf{Proposition}}[section]
\newtheorem{corollary}{\textbf{Corollary}}[section]
\newtheorem{remark}{\textbf{Remark}}[section]
\newtheorem{definition}{\textbf{Definition}}[section]

\allowdisplaybreaks[4]

\def\be{\begin{equation}}
\def\ee{\end{equation}}
\def\bea{\begin{eqnarray}}
\def\eea{\end{eqnarray}}
\def\bt{\begin{theorem}}
\def\et{\end{theorem}}
\def\bl{\begin{lemma}}
\def\el{\end{lemma}}
\def\br{\begin{remark}}
\def\er{\end{remark}}
\def\bc{\begin{corollary}}
\def\ec{\end{corollary}}
\def\bd{\begin{definition}}
\def\ed{\end{definition}}

\def\s0t{\sup _{0\leq \tau\leq t}}
\def\C0T{C([0,T];\,}

\def\DAS{D( A^{\frac{S}{2}})}
\def\DAS1{D( A^{ \frac{S+1}{2}})}

 \def\non{\nonumber }

\def \au {\rm}
\def \ti {\it}
\def \jou {\rm}
\def \bk {\it}
\def \no#1#2#3 {{\bf #1} (#3), #2.}
\def \eds#1#2#3 {#1, #2, #3.}
\begin{document}

\title{Long-time Behavior for a Nonlinear Plate Equation with Thermal Memory}

\author{
{\sc Hao Wu} \\ School of Mathematical Sciences, Fudan University \\
200433 Shanghai,
P.R.~China\\haowufd@yahoo.com}

\date{\today}

\maketitle

%%%%%%%%%%%%%%%%%%%%%%%%%%%%%%%%%%%%%%%%%%%%%%%%%%%%%%%%%%%%%%%%%%%%%%%%%%%%%%%%%%%%%%%%%%%%%%%%%%%%%%%%%%%%%%%%%%%%%%%%%%%

\begin{abstract}
\noindent We consider a nonlinear plate equation with thermal memory
effects due to non-Fourier heat flux laws. First we prove the
existence and uniqueness of global solutions as well as the
existence of a global attractor. Then we use a suitable \L
ojasiewicz--Simon type inequality to show the convergence of global
solutions to single steady states as time goes to infinity under the
assumption that the nonlinear term $f$ is real analytic. Moreover,
we
provide an estimate on the convergence rate.\\
\textbf{Keywords}: nonlinear plate equation, thermal memory, global attractor, convergence to equilibrium,
\L ojasiewicz--Simon inequality. \\
\textbf{Mathematics Subject Classification 2000}: 35B40, 35B41,
35B45.
\end{abstract}

%%%%%%%%%%%%%%%%%%%%%%%%%%%%%%%%%%%%%%%%%%%%%%%%%%%%%%%%%%%%%%%%%%%%%%%%%%%%
\section{Introduction}
\setcounter{equation}{0} In this paper, we consider the following
nonlinear plate equation with thermal memory effects due to
non-Fourier heat flux laws
 \begin{equation}
  \left\{\begin{array}{l} \theta_t -\Delta u_t +\int_0^\infty \kappa(s)[-\Delta \theta(t-s)]ds  =0, \\
    u_{tt} -\Delta u_t+ \Delta(\Delta u +\theta) +f(u) =0,
   \end{array}
   \right.\label{1a}
 \end{equation}
 for $(t,x)\in \mathbb{R}^+\times \Omega$, subject to the boundary conditions
\begin{equation}
        \left\{\begin{array}{l} u(t)=\Delta u(t)=0, \qquad t\geq 0,\ x\in  \Gamma, \\
    \theta(t)=0,\qquad \qquad \ \ t\in \mathbb{R},\ x\in  \Gamma,
   \end{array}
   \right.\label{1b}
 \end{equation}
 and initial conditions
\begin{equation}
 \begin{array}{l} u(0)=u_0,\  u_t(0)=v_0, \ \theta(0)=\theta_0,\quad x \in \Omega,\\
 \theta(-s)=\phi(s), \qquad (s, x) \in \mathbb{R}^+\times \Omega.
   \end{array}
   \label{1c}
\end{equation}
Here, $\Omega\in \mathbb{R}^2$ is a bounded domain with smooth
boundary $\Gamma$, $\theta$ represents the temperature variation
from the equilibrium reference value while $u$ is the vertical
displacement of the plate. Function
$\phi:\mathbb{R}^+\times\Omega\mapsto \mathbb{R}$ is called the
initial past history of temperature. The memory kernel
$\kappa:\mathbb{R}^+\mapsto\mathbb{R}$ is assumed to be a positive
bounded convex function vanishing at infinity. For the sake of
simplicity, we set all the physical constants to be one.

 Recently, evolution equations under various non-Fourier heat
 flux laws have attracted interests of many mathematicians
 (cf. \cite{AF,AP,GGP1,GP,CMP,FLR, GAA,GHS,GP01,GPV,GRP, Gra,Mola,Molaphd} and references cited therein).
  Let $\textbf{q}$ be the heat flux vector.
 According to the Gurtin--Pinkin theory \cite{GP1}, the
 linearized constitutive equation of $\textbf{q}$ is given by
 \begin{equation}
\textbf{q}(t) =-\int_0^\infty \kappa(s)\nabla \theta(t-s) ds,
\label{gp}
 \end{equation}
where $\kappa$ is the heat conductivity relaxation kernel. The
presence of convolution term in \eqref{gp} entails finite
propagation speed of thermal disturbances, so that in this case the
corresponding equation is of hyperbolic type. It is easy to see that
\eqref{gp} can be reduced to the classical Fourier law
$\textbf{q}=-\nabla \theta$ if $\kappa$ is the Dirac mass at zero.
Besides, if we take
\begin{equation}
\kappa(s)=\frac{1}{\sigma}e^{-\frac{s}{\sigma}},\qquad \sigma>0,
\end{equation}
and differentiate \eqref{gp} with respect to $t$, we can (formally)
arrive at the so called Cattaneo--Fourier law (cf.
\cite{He,Jo1,Jo2})
\begin{equation}
\sigma\textbf{q}_{t}(t)+\textbf{q}(t)=-\nabla\theta(t).
\end{equation}
On the other hand, evolution equations under Colemann--Gurtin theory
for the heat conduction (cf. \cite{CG}) have also been studied
extensively (see, for instance \cite{AF,AP,GPV,GGP}). There the heat
flux $\textbf{q}$ depends on both the past history and on the
instantaneous of the gradient of temperature:
\begin{equation}
\textbf{q}(t)= -K_I\nabla \theta(t) -\int_0^\infty \kappa(s)\nabla
\theta(t-s) ds,\label{CG}
\end{equation}
where $K_I>0$ is the instantaneous diffusivity coefficient.

 There are a lot of work on thermoelastic plate equations in
 the literature. For linear thermoelastic plate equations without memory effects in heat conduction,
 exponential stability of
 the associated  $C_0$-semigroups has been proven
 under different boundary conditions (cf. \cite[Section 2.5]{LZ}, \cite{RR1,RR}).
 On the other hand, when the heat flux is modeled
 by non-Fourier laws, wellposedness and stability for the corresponding linear thermoelastic plate equations
 have been investigated in several recent papers
 (cf. \cite {GP,GRP} and references cited therein).

 In this paper, we consider the nonlinear problem \eqref{1a}--\eqref{1c}. Asymptotic behavior of global solutions to
nonlinear isothermal plate equations has been considered before. We
may refer to \cite{HJ99,HR02}, where convergence to equilibrium as
$t\to \infty$ was obtained by the well known \L ojasiewicz--Simon
approach under the assumption that the nonlinearity is real
analytic. However, to the best of our knowledge, there are few
results on the long-time behavior of global solutions to nonlinear
plate equations with thermal memory like \eqref{1a}--\eqref{1c}.
This is just the main goal of the present paper. First, we prove the
existence and uniqueness of global solutions to
\eqref{1a}--\eqref{1c}. Then we derive some uniform estimates which
yields the precompactness of the solutions and furthermore the
existence of a global attractor. Finally, combining some techniques
for evolution equations with memory and for plate equations, we are
able to prove the convergence of global solutions to single steady
states as time goes to infinity via a suitable \L ojasiewicz--Simon
type inequality. Moreover, we obtain some estimates on convergence
rate. Further investigations concerning the infinite dimensional
system associated with our problem such as existence of exponential
attractors etc can be made by adapting the arguments in recent
papers \cite{GGP,Molaphd}.

 Our problem \eqref{1a}--\eqref{1c} is an evolution
 system with memory. It is well known in the literature that it would be more
 convenient to work in the history space setting by introducing a new variable $\eta$ called summed past history of $\theta$.
  This approach has been proven to be very effective
 in analyzing such kind of evolution systems (cf. \cite{GGP1,GGP,AF,AP,GAA,GPV,Gra,Mola,Molaphd,GP01,GRP}). On the other hand,
 it has been pointed out in the previous literature that when memory effects are present, the additional variable $\eta$ does not
 enjoy any regularizing effect. As a result, to
 ensure the precompactness of the trajectory, we have to make suitable decomposition of the solution which is typical for
 dissipative systems. To overcome the lack of compactness of the history space $\mathcal{M}$ in which the variable $\eta$ exists, an
 \textit{ad hoc} compactness lemma  will be used (cf. \cite{GGP1, Gra, GGP}).

 Comparing with the Colemann--Gurtin law (cf.
 \cite{AF,AP,GPV}), the dissipation in temperature $\theta$ for our system is only due to the memory effect, which is rather
 weak. The stronger dissipation provided in the Colemann--Gurtin law  would make the problem easier to be dealt with.
  For instance, we can refer to \cite{AF} in which the authors
 considered a nonisothermal phase--field system with \eqref{CG}
 and proved convergence to equilibrium for global solutions by the \L
 ojasiewicz--Simon approach (see also \cite{AP} for a conserved phase--field model). To overcome the difficulty due to such a weaker dissipation under the
 Gurtin--Pinkin law \eqref{gp}, it is necessary to introduce
 a suitable additional functional which may vary from problem to problem to
 produce some
 new dissipations (cf. \cite{GGP,GGP1,Gra,GAA, Mola} and the
 references cited therein). By using this idea, convergence to equilibrium for a nonisothermal Cahn--Hilliard
 equation was proven in \cite{Mola} and in \cite{GAA}
 a nonconserved phase--field model of Caginalp type consisting of two coupled integro-partial differential
equations was successfully treated.  Besides, in order to prove the
convergence result for our problem, we have to make use of an
extended \L ojasiewicz--Simon type
 inequality associated with a fourth order operator, which can be derived from the abstract result in \cite{HJ99}.
 Due to the structure of \eqref{1b}, the standard \L
ojasiewicz--Simon approach used in the parabolic case must be
modified by introducing an appropriate auxiliary functional (see
Section~5) which usually depends on the problem under consideration
(cf. \cite{HJ99,WGZ1,Molaphd,Mola,GAA} and references therein). In
our case, the required auxiliary functional is formed by adding two
perturbations to the original Lyapunov functional of system
\eqref{1a}--\eqref{1c} and coefficients of those perturbations
should be chosen properly. As far as the convergence rate is
concerned, it is known that an estimate in certain (lower order)
norm can usually be obtained directly from the \L ojasiewicz--Simon
approach (see, e.g., \cite{Z,HJ01}). Then one straightforward way to
get estimates in higher order norms is using interpolation
inequalities (cf. \cite{HJ01,Mola,GHS}) and, consequently, the decay
exponent deteriorates. We shall show that by using suitable energy
estimates and constructing proper differential inequalities, it is
possible to obtain the same estimates on convergence rate in both
higher and lower order norms. In particular, we find that as long as
uniform estimates in certain norm can be obtained, we are able to
prove convergence rate in the corresponding norm without loss in the
decay exponent. In our case, we can also avoid using the
decomposition argument used in \cite{Mola} for this purpose. This
technique has been successfully applied to other problems as well
(cf. \cite{W07,WGZ1,JJ1, WGZ3}) and it could be used to improve some
previous results in the literature (e.g., \cite{WGZ2,GHS, GAA,
Mola}). At last we show that actually better results on convergence
rate for problem \eqref{1a}--\eqref{1c} can be obtained if we use
the decomposition of the trajectory $z=z_D+z_C$ (see Section 4).
More precisely, the decay part $z_D$ converges to zero exponentially
fast while the compact part $z_C$ converges to equilibrium in a
higher order norm with the same rate as for the whole trajectory.

  The remaining part of this paper is organized as follows. In Section 2, we introduce the functional setting, the main results
  of this paper and some technical lemmas. Wellposedness of problem \eqref{1a}--\eqref{1c} is proven in Section 3. Section 4 is
  devoted to the uniform estimates and precompactness as well as the existence of a global attractor. In the final Section
  5, we prove the convergence of global solutions to
single steady states as time goes to infinity and obtain an estimate
on convergence rate.

%%%%%%%%%%%%%%%%%%%%%%%%%%%%%%%%%%%%%%%%%%%%%%%%%%%%%%%%%%%%%%%%%%%%%%%%%%%%%%%%%%%%%%%%%%%%%%%%%%%%%%%%%%%%%
\section{Preliminaries and Main Results}
\setcounter{equation}{0} We shall work under the functional settings
used in e.g. \cite{GRP}. Consider the positive operator $A$ on
$L^2(\Omega)$ defined by $A=-\Delta$ with domain
$D(A)=H^2(\Omega)\cap H^1_0(\Omega)$. Consequently, for $r\in
\mathbb{R}$ we can introduce the Hilbert spaces $V^r=D(A^{r/2})$,
endowed with the inner products
$$<w_1, w_2>_{V^r}=<A^{r/2}w_1,A^{r/2}w_2>,\qquad \forall\ w_1,w_2\in V^r, $$
where $<\cdot,\cdot>$ denotes the inner product in $L^2(\Omega)$. It
is easy to see that the embedding $ V^{r_1}\hookrightarrow V^{r_2}$
is compact for $r_1>r_2$. In what follows, we shall denote the norm
in $L^2(\Omega)$ by $\|\cdot\|$  for the sake of simplicity. \\
We suppose that $\kappa$ is vanishing at $\infty$. Moreover,
denoting
$$\mu(s)=-\kappa'(s),$$ we make the following assumptions on
$\mu$. \\
(H1) $\mu\in W^{1,1}(\mathbb{R}^+),$\\
(H2) $\mu(s)\geq 0, \qquad \mu'(s)\leq 0,\qquad \forall\ s\in
\mathbb{R}^+,$\\
(H3) $\mu'(s)+\delta\mu(s)\leq 0,\qquad \text{for some}\ \delta>0,\
\forall\ s\in \mathbb{R}^+,$\\
(H4) $\kappa(0)=\int_0^\infty\mu(s)ds:=\kappa_0> 0$.

 From recent
work \cite{CMP,pa,GAA} and the references cited therein, assumptions
on $\mu$ might be properly weakened and our results still hold. Our
results also hold under the assumptions made in \cite{Mola,Molaphd}
where (H1) (cf. \cite{GAA}) is replaced by $\mu\in
C^1(\mathbb{R}^+)\cap L^1(\mathbb{R}^+)$. In that case, $\mu$ is
allowed to be unbounded in a right neighborhood of 0 and this can be
handled by introducing a "cut-off" function near the origin.

For the nonlinear term $f$, we assume that\\
(F1) $f(s)\in C^2(\mathbb{R})$. \\
(F2)
\begin{equation}
\liminf_{|s|\rightarrow +\infty} \frac{f(s)}{s}>-\frac{1}{C_\Omega},
\non
\end{equation}
where $C_\Omega $ is the best constant depending only on $\Omega$
such that
$$\|w\|^2_{L^2(\Omega)}\leq C_\Omega\|Aw\|^2_{L^2(\Omega)}.$$
In order to prove the convergence to steady states, instead of (F1),
we assume \\
(F1)' $f(s)$ is real analytic in $s\in \mathbb{R}$.

We will also make use of the Poincar\'e inequality
$$ \|w\|\leq C_P\|\nabla w\|,\qquad w\in H_0^1(\Omega),$$
where $C_P$ is a positive constant depending only on $\Omega$.

In view of (H1)(H2), we introduce the weighted Hilbert spaces for
$r\in\mathbb{R}$,
$$\mathcal{M}^r=L^2_\mu(\mathbb{R}^+;V^r), $$
with inner products
$$<\eta_1,\eta_2>_{\mathcal{M}^r}=\int_0^\infty \mu(s)<A^{r/2}\eta_1(s),  A^{r/2}\eta_2(s)>ds.$$
Here we notice that the embeddings $\mathcal{M}^{r_1}\hookrightarrow
\mathcal{M}^{r_2}$, for $r_1> r_2$, are continuous but not compact
(cf. \cite{GGP,GGP1}).

Finally, we define the product Hilbert spaces
$$\mathcal{V}^r=V^{2+r}\times V^r\times V^r\times\mathcal{M}^{1+r},\qquad  r\in\mathbb{R},$$
with norm
$$
\|z\|^2_{\mathcal{V}^r}=\|A^{(2+r)/2}z_1\|^2+\|A^{r/2}z_2\|^2+\|A^{r/2}z_3\|^2+\|z_4\|^2_{\mathcal{M}^{1+r}},
$$
for all $z=(z_1,z_2,z_3,z_4)^T\in \mathcal{V}^r$.

It is convenient to work in the history space setting by introducing
the so-called summed past history of $\theta$ which is defined as
follows (cf. \cite{D, GP, GRP}),
 \be \eta^t(s)=\int_0^s\theta(t-y)dy,\qquad (t,s)\in [0,\infty) \times \mathbb
 {R}^+.
 \ee
The variable $\eta^t$ (formally) satisfies the linear equation
 \be
 \eta_t^t(s)+\eta_s^t(s)=\theta(t),\qquad \text{in}\ \Omega, \ (t,s)\in
 \mathbb{R}^+\times\mathbb{R}^+,\label{2.2}
 \ee
subject to the boundary and initial conditions
 \be \eta^t(0)=0,\qquad \text{in} \ \Omega, \ t\geq 0,\label{iniet}
 \ee
 \be
 \eta^0(s)=\eta_0(s)=\int_0^s\phi(y)dy, \qquad \text{in} \ \Omega, \ s\in \mathbb{R}^+.
 \ee
We introduce a linear operator $T$  on $\mathcal{M}^1$ defined by
 \be T\eta=-\eta_s, \qquad \eta\in D(T)
 \ee
 with domain
 \be
 D(T)=\left\{\eta\in \mathcal{M}^1\left|\eta_s\in \mathcal{M}^1,\eta(0)=0\right.\right\},
 \ee
here and in above $\eta_s$ is the distributional derivative of
$\eta$ with respect to internal variable $s$.

As in \cite{Gra,GAA}, we notice that an integration by parts in time
of the convolution products appearing in the equation for $\theta$
leads to
\begin{equation}
   \left\{\begin{array}{l} \theta_t -\Delta u_t -\int_0^\infty \mu(s)\Delta\eta^t(s)ds  =0, \\
    u_{tt} -\Delta u_t+ \Delta(\Delta u +\theta) +f(u) =0.
   \end{array}
   \right.\label{1a1}
 \end{equation}
Let us now introduce the vector
$$z(t)=\left(u(t),v(t),\theta(t),\eta^t\right)^T,$$
and denote the initial data by
 $$
 z_0=\left(u_0,v_0,\theta_0,\eta_0\right)^T\in \mathcal{V}^0.
 $$ Our problem \eqref{1a1}\eqref{1b}\eqref{1c} can be translated into the nonlinear abstract evolution
equation in $\mathcal{V}^0$,
\begin{equation}
  \left\{\begin{array}{l} z_t=Lz+G(z), \\
    z(0)=z_0,
   \end{array}
   \right.\label{evo}
 \end{equation}
with
\begin{equation}
G(z)=(0,-f(u),0,0)^T.
\end{equation}
Here the linear operator $L$ is defined as
\begin{equation}
L\left(\begin{array}{l} u \\v\\
\theta\\\eta\end{array}\right)=\left(\begin{array}{l} v
\\-Av-A(Au-\theta) \\ -Av-\int_0^\infty\mu(s)A\eta^t(s)ds\\ \theta+T\eta\end{array}\right),
\end{equation}
with domain
\begin{equation}
D(L)=\left\{z\in\mathcal{V}^0\left|\begin{array}{l} v, \ Au-\theta\in V^2 \\ \theta\in \mathcal{M}^1,\\
\int_0^\infty\mu(s)A\eta^t(s)ds  \in V^0\\ \eta\in
D(T)\end{array}\right.\right\}.
\end{equation}
 \br System \eqref{evo} is obtained through formal integration by
parts, however one can show that it is in fact equivalent to the
original problem \eqref{1a}--\eqref{1c} (cf. \cite{GRP}).
 \er

Now we are ready to state the main results of this paper
\begin{theorem}\label{MAIN1}
Let (H1)--(H4) and (F1)(F2) hold. The semigroup associated with
problem \eqref{evo} in $\mathcal{V}^0$ possesses a compact global
attractor $\mathcal {A}$ in $\mathcal{V}^0$.
\end{theorem}

\begin{theorem}\label{MAIN2}
Let (H1)--(H4) and (F1)'(F2) hold. Then for any $z_0=(u_0, v_0,
\theta_0,\eta_0)^T \in \mathcal{V}^0$, there exists $u_{\infty}$
being a solution to the following equation
\begin{equation}\left\{\begin{array}{l} A^2u_\infty+f(u_\infty)=0,\quad x\in \Omega, \\
    u_\infty=\Delta u_\infty =0,\quad x\in \Gamma,
   \end{array}
   \right.\label{sta}
\end{equation}
such that as $t \to \infty$
\begin{equation}\label{convu}
u(t) \to u_\infty \quad \text{in} \quad V^2,
\end{equation}
\begin{equation}\label{conv1}
v(t) \to 0, \quad \theta(t)\to 0 \qquad \text{in} \quad L^2(\Omega),
\end{equation}
\begin{equation}\label{conveta}
\eta^t \to 0 \quad \text{in} \quad \mathcal{M}^1.
\end{equation}
Moreover, there exists a positive constant $C$ depending on the
initial data such that
\begin{equation}\label{rateCSE}
\|u(t)
-u_\infty\|_{V^2}+\|v(t)\|+\|\theta(t)\|+\|\eta^t\|_{\mathcal{M}^1}
\leq C(1+t)^{-\frac{\rho}{(1-2\rho)}}, \quad \forall\ t \geq 0,
\end{equation}
with $\rho \in (0,1/2)$ being the same constant as in the \L
 ojasiewicz--Simon inequality (see Lemma \ref{ls}).
 \end{theorem}

\br With minor modifications, corresponding results can be proven
for equations under various other type of non-Fourier heat
conduction laws:
\begin{equation}
 \theta_t +c_1\theta-c_2\Delta\theta -\Delta u_t +\int_0^\infty\kappa(s)
 [c_3\theta(t-s)-\Delta \theta(t-s)]ds  =0,\label{onf}
\end{equation}
with $c_1,c_2,c_3$ being nonnegative constants. When $c_1>0$, we
have a (dissipative) term $c_1\theta$ in \eqref{onf}, which is
arising from the assumption that besides the heat flux, the thermal
power depends on the past history of $\theta$ (cf. \cite{GP}). The
case $c_2>0$ corresponds to the Colemann--Gurtin theory as mentioned
before. Moreover, we may refer to \cite{GP,L} for the case $c_3>0$.
Although there might be additional terms like $c_1\theta$,
$-c_2\Delta \theta$ and $\int_0^\infty\kappa(s)
 c_3\theta(t-s) ds$ in the equation, these terms
provide stronger dissipations on $\theta$ from the mathematical
point of view, which make the extensions of our results possible.\er

For reader's convenience we report below some helpful technical
lemmas which will be used in this paper. The first one is a
frequently used compactness lemma for the spaces $\mathcal{M}^r$
(cf. \cite[Lemma 2.1]{GRP}). \bl\label{com4} Let
$\mathbb{T}_{\eta}(y)$ be defined as follows
\begin{equation}
\mathbb{T}_{\eta}(y)=\int_{(0,1/y)\cup(y,\infty)}\mu(s)\|A^{1/2}\eta(s)\|^2ds,\qquad
y\geq 1.
\end{equation}
If $\mathcal{C}\subset \mathcal{M}^1$ satisfiy \\(i) $\sup_{\eta\in
\mathcal{C}}\|\eta\|_{\mathcal{M}^2}<\infty,$ \\(ii) $\sup_{\eta\in
\mathcal{C}}\|T\eta\|_{\mathcal{M}^1}<\infty,$\\ (iii)
$\lim_{y\rightarrow \infty}\left(\sup_{\eta\in
\mathcal{C}}\mathbb{T}_{\eta}(y)\right)=0,$ \\ then $\mathcal{C}$ is
relatively compact in $\mathcal{M}^1$. \el The following lemma can
be found in \cite{BP}.
 \bl\label{tl1} Let $X$ be a Banach space and $Z\in C([0,\infty),X)$. Let $E:X\rightarrow \mathbb{R}$ be a
function bounded from below such that $E(Z(0))\leq M$ for $Z\in X$.
If
$$\frac{d} {dt}E(Z(t))+\delta \|Z(t)\|_X^2\leq k,$$
 for some $\delta\geq 0$ and $k\geq 0$ independent of $Z$, then for all $\varepsilon>0$ there is
 $t_0=t_0(M,\varepsilon)>0$ such that
 $$      E(Z(t))\leq \sup_{\xi\in X}\left\{E(\xi):\delta\|\xi\|_X^2\leq k+\varepsilon\right\},\qquad \forall\ t\geq t_0.$$
 \el

%%%%%%%%%%%%%%%%%%%%%%%%%%%%%%%%%%%%%%%%%%%%%%%%%%%%%%%%%%%%%%%%%%%%%%%%%%%%%%%%%%%%%%%%%%%%%%%%%%%%%%%%%%%%%
\section{Wellposedness}
\setcounter{equation}{0} By using the semigroup approach, we are
able to prove the existence and uniqueness of global solution to
system \eqref{evo}.
\begin{theorem}\label{segr} Suppose that assumptions (H1)(H2) and (F1)(F2) hold.
Then for any initial data $z_0=(u_0,v_0,\theta_0,\eta_0)^T\in
\mathcal{V}^0$, system \eqref{evo} admits a unique global solution
$z(t)\in C([0,+\infty), \mathcal{V}^0)$.
\end{theorem}
\begin{proof}
We apply the semigroup theory (see e.g., \cite[Theorem 2.5.4, Theorem 2.5.5]{Z}). \\
Since
\begin{equation}
 <T\eta,\eta>_{\mathcal{M}^1}=\frac12\int_0^\infty\mu'(s)\|A^{1/2}\eta(s)\|^2ds\leq 0,\qquad \forall\  \eta\in D(T),\label{de1}
\end{equation}
it's easy to see that
\begin{equation}
<Lz,z>_{\mathcal{V}^0}\ =\ -\|\nabla
v\|^2+<T\eta,\eta>_{\mathcal{M}^1}\leq 0,\qquad \forall \ z\in
D(L).\label{de2}
\end{equation}
By a similar argument in \cite[Section 3]{GP} (see also \cite{GRP})
we can show that $I-L: D(L) \mapsto \mathcal{V}^0$ is onto. Thus $L$
is a m-accretive operator. On the other hand, by the Sobolev
embedding Theorem, for any $z_1, z_2\in \mathcal{V}^0$ with $
\|z_1\|_{\mathcal{V}^0}\leq M$, $\ \|z_2\|_{\mathcal{V}^0}\leq M$,
there exists a constant $L_M>0$ depending on $M$ such that
$$\|G(z_1)-G(z_2)\|_{\mathcal{V}^0}\leq  L_M\|u_1-u_2\|_{V^2}\leq L_M\|z_1-z_2\|_{\mathcal{V}^0}.$$
 Therefore, $ G(z)$ is a nonlinear operator form $\mathcal{V}^0$ to $\mathcal{V}^0$ satisfying the
 local Lipschitz condition. Consequently, local existence of a unique mild solution $z(t)\in C([0,T], \mathcal{V}^0)$
 follows from \cite[Theorem 2.5.4]{Z}.

Next we prove the global existence. Taking inner product of
\eqref{evo} and $z$ in $ \mathcal{V}^0$, we get
\begin{equation}
\frac{d}{dt}\left(\frac{1}{2}\|z(t)\|^2_{\mathcal{V}^0}+\int_\Omega
F(u)dx\right) +\|\nabla
v\|^2-\frac12\int_0^\infty\mu'(s)\|A^{1/2}\eta^t(s)\|^2ds=
0,\label{ly}
\end{equation}
where $F(u)=\int_0^u f(y)dy$.\\
 Assumption (F2) implies that there exist constants $\delta\in (0,1) $ and
    $N=N(\delta)>0$ such that (cf. \cite{CEL})
   \begin{equation} F(s)\geq -\frac{1-\delta}{2C_\Omega}
    s^2, \qquad \text{for} \ |s|\geq N. \nonumber
    \end{equation}
    To see this, let $M$ be a positive
    constant such that $f(z)/z + \frac{1}{C_\Omega} \geq \frac{2\delta}{C_\Omega}$ for $|z|\geq M$
    and certain $\delta\in( 0, 1)$. Then we have
     \bea
    F(s)+\frac{1}{2C_\Omega}s^2 & =& \int_0^M\left(\frac{f(z)}{z}+ \frac{1}{C_\Omega}\right) z
    dz + \int_M^s\left(\frac{f(z)}{z}+ \frac{1}{C_\Omega}\right) z dz \nonumber\\
    & \geq & C+\frac{2\delta}{C_\Omega}\left(\frac{s^2}{2}-\frac{M^2}{2}\right) \geq
    \frac{\delta}{2C_\Omega}s^2
     \eea
     for $$s^2\geq \max\left\{2M^2\!-\frac{2C_\Omega C}{\delta},\  0\right\}:=N^2.$$ For
    negative $s$ one can repeat the same computation with $M$ replaced
    by $-M$. \\
    Now we have \be \int_\Omega F(u) dx= \int_{|u|\leq N}
    F(u)dx + \int_{|u|> N}F(u)dx \geq
    -\frac{1-\delta}{2C_\Omega}\int_\Omega u^2dx + C(|\Omega|,f)\ee
    where $C(|\Omega|,f)=|\Omega|\displaystyle{\min_{|s|\leq N}}F(s)$.\\
    By the definition of $C_\Omega$ in (F2) we can deduce
     \begin{equation} \int_\Omega F(u)dx \geq  - \frac{1-\delta}{2} \|A u\|^2 +
     C(|\Omega|,f).
          \end{equation}
   This implies that for any $\epsilon \in(0,\delta]$ there holds
   \begin{eqnarray}
   \frac{1}{2}\|z(t)\|^2_{\mathcal{V}^0}+\int_\Omega F(u)dx &=&  \frac{\epsilon}{2}\|z(t)\|^2_{\mathcal{V}^0}+
   \frac{1-\epsilon}{2}\|z(t)\|^2_{\mathcal{V}^0}+\int_\Omega F(u)dx\nonumber\\
    &\geq& \frac{\epsilon}{2}\|z(t)\|^2_{\mathcal{V}^0}+ C(|\Omega|,f).\label{ee}
   \end{eqnarray}
  As a result,
  \begin{equation}  \frac{1}{2}\|z(t)\|^2_{\mathcal{V}^0}\leq \frac{1}{\epsilon}
  \left(\frac{1}{2}\|z(t)\|^2_{\mathcal{V}^0}+\int_\Omega F(u)dx- C(|\Omega|,f)\right).\label{e}
  \end{equation}
     Integrating \eqref{ly} with respect to $t$,  we infer from
     \eqref{e} that
\begin{equation}
   \|z(t)\|^2_{\mathcal{V}^0}\leq C(\|z_0\|_{\mathcal{V}^0}, |\Omega|, f),\qquad \forall \ t\geq 0.\label{w1}
\end{equation}
This uniform estimate together with \cite[Theorem 2.5.5]{Z} yields
the global existence, i.e., $z(t)\in C([0,+\infty), \mathcal{V}^0)$.
Moreover, it is not difficult to check that for any $z_{01},
z_{02}\in \mathcal{V}^0$, the corresponding global solutions
$z_1(t), z_2(t)$ satisfy
\begin{equation}
\| z_1(t)- z_2(t)\|^2_{\mathcal{V}^0}\leq C_T \| z_{01}-
z_{02}\|^2_{\mathcal{V}^0},\qquad 0\leq t\leq T,\label{uni}
\end{equation}
for all $T\geq 0$, where $C_T$ is a constant depending on the norms of $ z_{01}, z_{02}$ in $\mathcal{V}^0$ and $T$.\\
The proof is complete.
\end{proof}
\br
 From the above theorem, we can see that the solution to our problem \eqref{evo} defines a
strongly continuous semigroup $S(t)$ on the phase space
$\mathcal{V}^0$ such that $S(t)z_0=z(t)$. \er

%%%%%%%%%%%%%%%%%%%%%%%%%%%%%%%%%%%%%%%%%%%%%%%%%%%%%%%%%%%%%%%%%%%%%%%%%%%%%%%%%%%%%%%%%%%%%%%%%%
\section{Precompactness of Trajectories and Global Attractor}
\setcounter{equation}{0} In this section, we will first prove (i)
uniform estimate of the solution which also indicates the existence
of an absorbing set, (ii) precompactness of trajectory $z(t)$. In
what follows, we shall exploit some formal \textit{a priori}
estimates which can be justified rigorously by the standard density
argument.
\begin{lemma}
\label{abs1} Let (H1)--(H4) and (F1)(F2) hold. There exists a
positive constant $R_0$ such that the ball
$$\mathcal{B}_0:=\{z\in \mathcal{V}^0\ |\ \| z\|_{\mathcal{V}^0}\leq R_0\}$$
 is an absorbing set. Namely, for any bounded set $\mathcal{B}\in \mathcal{V}^{0}$,
 there is $t_0=t_0(\mathcal{B})\geq 0$ such that
  $S(t)\mathcal{B}\subset \mathcal{B}_0$ for every $t\geq t_0$.
\end{lemma}
\begin{proof}
Multiplying the second equation in \eqref{1a} by $\varepsilon^2 u$,
integrating on $\Omega$ and adding the result to \eqref{ly}, we get
\begin{eqnarray}
&& \frac{d}{dt}\left(
\frac{1}{2}\|z(t)\|^2_{\mathcal{V}^0}+\frac{\varepsilon^2}2\|\nabla
u\|^2+\int_\Omega F(u)dx+ \varepsilon^2
\int_\Omega uvdx\right)+\|\nabla v\|^2\non\\
&&\ \ -\varepsilon^2\|v\|^2-\frac12\int_0^\infty\mu'(s)\|A^{1/2}\eta^t(s)\|^2ds+ \varepsilon^2\|Au\|^2\nonumber\\
&=&  -\varepsilon^2\int_\Omega f(u)udx+\varepsilon^2\int_\Omega
\theta A udx.\label{ab1}
\end{eqnarray}
In order to apply Lemma \ref{tl1}, we need more dissipation on the
left hand side of \eqref{ab1}. To such an aim, we introduce the
following functional (cf. \cite{Gra,Mola,GAA} and the references
cited  therein)
\begin{equation}
J(t):=-\int_0^\infty \mu(s)<\theta(t),\eta^t(s)>ds.
\end{equation}
It turns out from the H\"older inequality and (H1) that
\begin{eqnarray}
|J(t)|&\leq& \|\theta(t)\|\int_0^\infty
\mu(s)<\eta^t(s),\eta^t(s)>^\frac12ds\non\\
 & \leq&
\frac12\|\theta(t)\|^2+\frac12\|\eta^t(s)\|_{\mathcal{M}^0}^2\int_0^\infty\mu(s)ds\non\\
& \leq& C\|z(t)\|^2_{\mathcal{V}^0}.\label{JJ}
\end{eqnarray}
Besides, a direct calculation yields (cf. \eqref{2.2}\eqref{1a1})
\begin{eqnarray}
\frac{d}{dt}J(t)&=& - \int_0^\infty
\mu(s)<\theta_t(t),\eta^t(s)>ds-\int_0^\infty
\mu(s)<\theta(t),\eta_t^t(s)>ds\non\\
&=& -\int_0^\infty\mu(s)<\Delta
u_t(t),\eta^t(s)>ds-\left\|\int_0^\infty\mu(s)
A^{1/2}\eta^t(s)ds\right\|^2-\kappa_0\|\theta\|^2\non\\
&& +\int_0^\infty\mu(s)<\theta(t),\eta^t_s(s)>ds.\label{4.4}
\end{eqnarray}
Terms on the right hand side of \eqref{4.4} can be controlled in the
following way:
\begin{eqnarray}
\left|-\int_0^\infty\mu(s)<\Delta u_t(t),\eta^t(s)>ds\right|&=&
\left|\int_0^\infty \mu(s)<\nabla v(t), \nabla
\eta^t(s)>ds\right|\non\\
&\leq&  \frac12\|\nabla v\|^2+
\frac{\kappa_0}{2}\|\eta^t\|_{\mathcal{M}^1}^2,\label{j1}
\end{eqnarray}
\begin{equation}
\left\|\int_0^\infty\mu(s) A^{1/2}\eta^t(s)ds\right\|^2\leq
\int_0^\infty \mu(s)ds \int_0^\infty
\mu(s)<A^{1/2}\eta^t,A^{1/2}\eta^t>ds\leq
\kappa_0\|\eta^t\|_{\mathcal{M}^1}^2,
\end{equation}
\begin{eqnarray}
\left|\int_0^\infty\mu(s)<\theta(t),\eta^t_s(s)>ds\right|&=&\left|
-\int_0^\infty
\mu'(s) <\theta(t),\eta^t(s)>ds\right|\non\\
&\leq& -\int_0^\infty \mu'(s)\|\theta(t)\|\|\eta^t(s)\|ds\non\\
&\leq& \frac{\kappa_0}{2}\|\theta\|^2- C_1 \int_0^\infty \mu'(s)
\|A^{1/2}\eta^t(s)\|^2ds,\label{j2}
\end{eqnarray}
where in \eqref{j2} we use (H1) that $\mu'$ is integrable (this can
be weakened as mentioned in the previous section) and $C_1>0$
depends on $\kappa_0$.
 Now we can conclude
\begin{equation}
\frac{d}{dt}J(t)+\frac{\kappa_0}{2}\|\theta\|^2 \leq \frac12\|\nabla
v\|^2+ C_2\|\eta^t\|_{\mathcal{M}^1}^2-C_1 \int_0^\infty \mu'(s)
\|A^{1/2}\eta^t\|^2ds,\label{aa2}
\end{equation}
where $C_2=\frac{3\kappa_0}{2}>0$. Multiplying \eqref{aa2} by
$2\varepsilon$ and adding it to \eqref{ab1} we obtain
\begin{eqnarray}
&& \frac{d}{dt}\left(
\frac{1}{2}\|z(t)\|^2_{\mathcal{V}^0}+\frac{\varepsilon^2}{2}\|\nabla
u\|^2+\int_\Omega F(u)dx+ \varepsilon^2\int_\Omega
uvdx+2\varepsilon J(t)\right)-\varepsilon^2\|v\|^2\non\\
& & \ +(1-\varepsilon)\|\nabla v\|^2 +\varepsilon\kappa_0
\|\theta\|^2-\left(\frac12-2C_1\varepsilon\right)\int_0^\infty\mu'(s)\|A^{1/2}\eta^t(s)\|^2ds+
\varepsilon^2\|Au\|^2\nonumber\\
&\leq&  -\varepsilon^2\int_\Omega f(u)udx+\varepsilon^2\int_\Omega
\theta A udx+ 2C_2\varepsilon
\|\eta^t\|_{\mathcal{M}^1}^2.\label{ab1a}
\end{eqnarray}
Define
$$ \Psi(z(t))=\frac12\|z(t)\|^2_{\mathcal{V}^0}+\frac{\varepsilon^2}{2}\|\nabla u(t)\|^2
+\int_\Omega F(u(t))dx+ \varepsilon^2\int_\Omega u(t)v(t)dx +
2\varepsilon J(t),$$ and
$$\mathcal{F}(\|u\|_{V^2})=|\Omega|\max_{|y|\leq \|u\|_{V^2}}|F(y)|.      $$
Due to (F1) and the Sobolev embedding Theorem $V^2\hookrightarrow
L^{\infty}(\Omega)$, we can see that $\mathcal{F}(s)$ is bounded for
$|s|\leq M,\ \forall M>0$. It follows from this fact and \eqref{ee}
that for all $z\in \mathcal{V}^0$ and $\varepsilon$ sufficiently
small there holds
\begin{equation}
C_3\|z\|_{\mathcal{V}^0}^2-C_4 \leq \Psi\leq
\|z\|_{\mathcal{V}^0}^2+\mathcal{F}(\|u\|_{V^2}),\label{phi}
\end{equation}
where $C_3$, $C_4$ are positive constants independent of $z$. \\
(F2) implies that there exist constants $\sigma > 0 $  and
$N=N(\sigma)>0$ satisfying
$$f(s)s\geq -\frac{1-\sigma}{C_\Omega}s^2,\qquad \forall\ |s|\geq N.$$
As a result,
\begin{equation}
\int_\Omega f(u)udx \geq -\frac{1-\sigma}{C_\Omega}\|u\|^2+ C(|\Omega|,f)\geq -(1-\sigma)\|Au\|^2
+ C(|\Omega|,f), \label{ab2aa}
\end{equation}
where $C(|\Omega|,f)=|\Omega|\displaystyle{\min_{|s|\leq N}}f(s)s$.\\
Moreover, from (H3), we can see that
\begin{equation}
 -\int_0^\infty\mu'(s)\|A^{1/2}\eta^t(s)\|^2ds \geq \delta \int_0^\infty\mu(s)\|A^{1/2}\eta^t(s)\|^2ds, \label{ab2}
\end{equation}
and by the H\"older inequality we have
\begin{equation}
\left|\varepsilon^2\int_\Omega \theta A u dx\right|\leq
\varepsilon\frac{\kappa_0}{2}\|\theta\|^2
+C_5\varepsilon^3\|Au\|^2.\label{ab2a}
\end{equation}
In \eqref{ab1a} and \eqref{ab2a} we take $\varepsilon$ small enough
such that
$$
0<\varepsilon\leq \min\left\{\frac14, \frac{1}{16C_1},
\frac{\delta}{16C_2}, \frac{\sigma}{2C_5}, \frac{1}{2C_P} \right\}.
$$
Then it follows from \eqref{ab1a}-\eqref{ab2a} that
\begin{equation}
 \frac{d}{dt}\Psi(t)+ \delta_0\|z(t)\|_{\mathcal{V}^0}^2\leq k,\label{ab3}
\end{equation}
where $\delta_0$ is a constant depending on $\varepsilon, \sigma,
\delta,\kappa_0$ and $k$ is a certain positive constant depending on
$\varepsilon, |\Omega|, f, \delta, \sigma, N(\sigma) $.

We infer from Lemma \ref{tl1} that there is $t_0=t_0(\mathcal{B})>0
$ such that
$$\Psi(z(t))\leq \sup_{\xi\in \mathcal{V}^0}\left\{\Psi(\xi): \delta_0\|\xi\|_{\mathcal{V}^0}^2\leq 1
+k\right\},\qquad \forall\ t\geq t_0,$$ which together with
\eqref{phi} implies the existence of absorbing set.

The proof is complete.
\end{proof}
Next we prove the precompactness of solutions to problem
\eqref{evo}. Since our system \eqref{evo} does not enjoy smooth
property as parabolic equations, it suffices to show that the
semigroup is asymptotically smooth (cf. \cite{T}). To accomplish
this, we make a decomposition of the flow into a uniformly stable
part and a compact part (cf. \cite{GGP,GGP1,GAA,Gra,Molaphd}).
Namely, we decompose the solution to \eqref{evo} with initial data
$z(0)=z_0\in \mathcal{V}^0$ as
$$z(t)=z_D(t)+z_C(t),$$ where
$z_D(t)=\left(u_D(t),v_D(t),\theta_D(t),\eta_D^t\right)^T$ and
$z_C(t)=\left(u_C(t),v_C(t) ,\theta_C(t),\eta_C^t\right)^T$ satisfy
\begin{equation}
  \left\{\begin{array}{l} \displaystyle \frac{d}{dt}z_D=Lz_D, \\
    z_D(0)=z_0,
   \end{array}
   \right.\label{evoD}
 \end{equation}
and
 \begin{equation}
  \left\{\begin{array}{l} \displaystyle \frac{d}{dt}z_C=Lz_C+G(z), \\
    z_C(0)=0.
   \end{array}
   \right.\label{evoC}
 \end{equation}
Similar to Theorem \ref{segr}, it is easy to check that system
\eqref{evoD} admits a unique mild solution $z_D(t)\in C([0,+\infty),
\mathcal{V}^0)$. Moreover, we have

\bl\label{zd} There exist constants $C,\delta_1>0$ such that the
solution $z_D$ of \eqref{evoD} fulfills
 \begin{equation}
 \|z_D(t)\|_{\mathcal{V}^0}\leq Ce^{-\frac{\delta_1}{2} t},\qquad \forall\  t\geq 0,\label{conD}
 \end{equation}
 where $C>0$ is a constant depending on $\|z_0\|_{\mathcal{V}^0}$.
\el
\begin{proof}
Let $E_D(z):\mathcal{V}^0\mapsto \mathbb {R}$ be defined as follows
\begin{equation}
E_D(z_D)=\|z_D(t)\|^2_{\mathcal{V}^0}+ \varepsilon^2 \|\nabla
u_D\|^2+ 2\varepsilon^2\int_\Omega u_Dv_Ddx+4\varepsilon J_D(t)
\end{equation}
with
\begin{equation}
J_D(t):=-\int_0^\infty \mu(s)<\theta_D(t),\eta_D^t(s)>ds.
\end{equation}
It is easy to see that for $\varepsilon>0$ sufficiently small
 \be
\frac12\|z_D\|^2_{\mathcal{V}^0}\leq E_D(z_D)\leq
2\|z_D\|^2_{\mathcal{V}^0}.\label{EZD} \ee
 Similar to the proof of Lemma \ref{abs1},
we can show that there exists $\delta_1>0$ such that
\begin{eqnarray}
&& \frac{d}{dt}E_D(z_D)+ 2\delta_1\|z_D\|^2_{\mathcal{V}^0}\leq
0.\label{ab5}
\end{eqnarray}
 As a consequence of \eqref{EZD}, we have
\begin{eqnarray}
&& \frac{d}{dt}E_D(z_D)+ \delta_1 E_D(z_D)\leq 0\label{aba5}
\end{eqnarray}
which yields
 \be E_D(z_D(t))\leq E_D(z_D(0))e^{-\delta_1t}, \quad \forall\ t\geq
 0.\label{EED}\ee
\eqref{conD} follows immediately from \eqref{EZD} and \eqref{EED}.
\end{proof}

Next we analyze $z_C$. For initial data $z_0\in \mathcal{V}^0$, we
can see that $z_C(t)=z(t)-z_D(t)$ belongs to a bounded set in
$\mathcal{V}^0$ for $t\geq 0$. In what follows we will show that
$z_C$ is more regular and actually it is uniformly bounded in
$\mathcal{V}^1$. \bl\label{zc}
 For all $z_0\in \mathcal{V}^0$, there exists $C> 0$ depending on $\| z_0\|_{\mathcal{V}^0}$ such that
\begin{equation}
 \|z_C(t)\|_{\mathcal{V}^1}\leq C,\qquad \forall\ t\geq 0.\label{zcb}
\end{equation}
\el
\begin{proof}
Taking the inner product of \eqref{evoC} and $Az_C$ in $
\mathcal{V}^0$, we have
 \bea
&&
\frac{d}{dt}\left(\frac12\|A^{1/2}z_C\|^2_{\mathcal{V}^0}+\int_\Omega
f(u)Au_Cdx\right)+\|A
v_C\|^2-\frac12\int_0^\infty\mu'(s)\|A\eta_C^t(s)\|^2ds\non\\
&=&\int_\Omega f'(u)vAu_C dx.\label{COM1}
 \eea
Multiplying the second equation in \eqref{evoC} by $Au_C$ and
integrating on $\Omega$ we get
 \bea
&& \frac{d}{dt}\left(\frac{1}{2}\|Au_C\|^2 +\int_\Omega v_CAu_Cdx
\right)+\|A^{3/2}u_C\|^2-\|A^{1/2} v_C\|^2 \non\\
&=& - \int_\Omega f(u)Au_Cdx+\int_\Omega
A\theta_CAu_Cdx.\label{COM2}
 \eea
 Let
\begin{equation}
J_C(t):=-\int_0^\infty
\mu(s)<A^{1/2}\theta_C(t),A^{1/2}\eta_C^t(s)>ds.
\end{equation}
Then in analogy to the argument in Lemma \ref{abs1}, we have
 \be
 |J_C(t)|\leq \frac12\|A^{1/2}\theta_C(t)\|^2+\frac12\|\eta_C^t(s)\|_{\mathcal{M}^1}^2\int_0^\infty\mu(s)ds
\leq C\|z_C(t)\|^2_{\mathcal{V}^1},\label{JC}
 \ee
and
\begin{equation}
\frac{d}{dt}J_C(t)+\frac{\kappa_0}{2}\|A^{1/2}\theta_C\|^2 \leq
\frac12\|A v_C\|^2+ C_6\|\eta_C^t\|_{\mathcal{M}^2}^2-C_7
\int_0^\infty \mu'(s) \|A\eta_C^t\|^2ds,\label{COM3}
\end{equation}
here $C_6, C_7>0$ depend on $\kappa_0$. Introduce the functional
$$\Phi(t)=\|A^{1/2}z_C(t)\|^2_{\mathcal{V}^0}+ 2\int_\Omega
f(u)Au_Cdx+ 2\varepsilon^2\int_\Omega v_CAu_Cdx + \varepsilon^2\|A
u_C\|^2+4\varepsilon J_C(t)+k,
$$
where $k\geq 0$ denotes a generic constant depending on $\|z_0\|_{\mathcal{V}^0}$.\\
It is easy to see that, if the constant $k$ appearing in the
definition of $\Phi$ is large enough and $\varepsilon$ is small
enough, there holds
\begin{equation}
\frac{1}{2}\|A^{1/2}z_C(t)\|^2_{\mathcal{V}^0}\leq \Phi\leq
2\|A^{1/2}z_C(t)\|^2_{\mathcal{V}^0}+C(\|z_0\|_{\mathcal{V}^0})+k.\label{bPhi}
\end{equation}
It follows from \eqref{COM1}\eqref{COM2}\eqref{COM3} that
 \bea
 &&\frac12\frac{d}{dt}\Phi(t)+ (1-\varepsilon)\|A
v_C\|^2-\left(\frac12-2C_7\varepsilon\right)\int_0^\infty\mu'(s)\|A\eta_C^t(s)\|^2ds
\non\\
&&-\varepsilon^2\|A^{1/2}v_C\|^2+\varepsilon^2\|A^{3/2}u_C\|^2
+\varepsilon
\kappa_0\|A^{1/2}\theta_C\|^2\non\\
&\leq& \int_\Omega f'(u)vAu_Cdx
+2C_6\varepsilon\|\eta_C^t\|_{\mathcal{M}^2}^2-\varepsilon^2
\int_\Omega f(u)Au_Cdx+\varepsilon^2\int_\Omega
A\theta_CAu_Cdx\non\\
&:=& I_1+I_2+I_3+I_4.\label{c1}
 \eea
The right hand side of \eqref{c1} can be estimated as follows. By
(F2) and the Sobolev embedding Theorem we get
\begin{eqnarray}
|I_1|&=& \left| \int_\Omega f'(u)vAu_Cdx\right| \leq  \|f'(u)\|_{L^\infty}\|v\|\|Au_C\|\nonumber\\
&\leq &C (|\Omega|,\|Au\|)\|v\|\|\|A^{3/2}u_C\|\nonumber\\
&\leq & \frac{\varepsilon^2}{4}\|A^{3/2}u_C\|^2+C.\label{ff3}
\end{eqnarray}
Besides,
 \be|I_3|= \left| -\varepsilon^2 \int_\Omega f(u)Au_Cdx\right|\leq
 \varepsilon^2 \|f(u)\|_{L^\infty}\|Au_C\|\leq
 \frac{\varepsilon^2}{4}\|A^{3/2}u_C\|^2+C,
 \ee
\begin{equation}
 |I_4|=\left| \varepsilon^2\int_\Omega Au_C A\theta_C dx\right|\leq
\frac{\varepsilon\kappa_0}{2}\|A^{1/2}\theta_C\|^2+C_8\varepsilon^3\|A^{3/2}u_C\|^2,\label{c2}
\end{equation}
where in the above three estimates, $C>0$ is a constant depending on
$\|z_0\|_{\mathcal{V}^0}$ and $C_8=\frac{1}{2\kappa_0}$. \\
In \eqref{c1}--\eqref{c2}, we take
$\varepsilon$ small enough such that
\begin{equation}
0< \varepsilon\leq \min\left\{\frac14,
\frac{1}{2C_P},\frac{\delta}{16C_6}, \frac{1}{16C_7},
\frac{1}{4C_8}\right\}.
\end{equation}
Consequently, we can obtain the following inequality
\begin{equation}
\frac{d}{dt}\Phi(t)+ \delta_2 \Phi(t) \leq C,\label{com}
\end{equation}
where $\delta_2>0$, $ C\geq 0$ are constants depending on $\varepsilon$ and $\|z_0\|_{\mathcal{V}^0}$.\\
\eqref{bPhi} and \eqref{com} yield
$$\|A^{1/2}z_C(t)\|_{\mathcal{V}^0}\leq C,\quad \forall\ t\geq 0.$$
The proof is complete.
\end{proof}
In order to obtain the required compactness, we have to take care of
the fourth component $\eta^t$. Embedding $
\mathcal{V}^1\hookrightarrow \mathcal{V}^0 $ is not compact because
embedding $  \mathcal{M}^2\hookrightarrow \mathcal{M}^1$ is not
compact in general. However, we have the following lemma whose proof
is becoming standard (cf. \cite{GGP1,GP01,Gra} and the references
therein). For the sake of completeness, we give a sketch of the
proof.
 \bl\label{com7} Let $\mathcal{C}=\bigcup_{t\geq
0}\eta_C^t$. Then $ \mathcal{C}$ is relatively compact in
$\mathcal{M}^1$. \el
\begin{proof}
It is obvious that $\mathcal{C} \subset \mathcal{M}^1$. According to Lemma  \ref{com4}, we need to verify
\begin{equation}\|\eta_C^t\|_{\mathcal{M}^2}\leq C, \qquad t\geq 0, \label{c4}
\end{equation}
\begin{equation}
\|T\eta_C^t\|_{\mathcal{M}^1}\leq C,\qquad t\geq 0. \label{c5}
\end{equation}
\begin{equation}
 \lim_{y\rightarrow \infty}\left(\sup_{t\geq 0 }\mathbb{T}_{\eta_C^t}(y)\right)=0. \label{c7}
\end{equation}
 \eqref{c4} follows from Lemma \ref{zc} immediately. Since $z_C(0)=0$, $\eta_C$ has the following explicit representation formula
  \begin{equation}
  \eta_C^t(s)=\left\{\begin{array}{ll} \int_0^s\theta_C(t-y)dy,& 0<s\leq t,  \\
   \int_0^t\theta_C(t-y)dy, & s>t.
   \end{array}
   \right.\label{ec}
 \end{equation}
Differentiating it with respect to $s$ yields
\begin{equation}
  T\eta_C^t(s)=\left\{\begin{array}{ll} -\theta_C(t-s),& 0<s\leq t,  \\
   0, & s>t.
   \end{array}
   \right.
 \end{equation}
Thanks to (H4) and Lemma \ref{zc}, we have
\begin{equation}
\int_0^\infty\mu(s)\|A^{1/2}T\eta_C^t(s)\|^2ds=\int_0^t\mu(s)\|A^{1/2}T\eta_C^t(s)\|^2ds=
\int_0^t\mu(s)\|A^{1/2}\theta_C(t-s)\|^2ds\leq C,
\end{equation}
which yields \eqref{c5}. This also implies that $\eta^t_C\in H^1_\mu(\mathbb{R}^+;V^1)$.\\
\eqref{ec} and Lemma \ref{zc} imply that
$$\|A^{1/2}\eta_C^t(s)\|^2\leq C(1+s^2), \qquad \forall\ s>0.$$
For $y\geq 1$, we define
$$I(y)=C\int_{(0,1/y)\cup(y,\infty)}\mu(s)(1+s^2)ds.$$
It's obvious that (see the definition of $\mathbb{T}$)
\begin{equation}\mathbb{T}_{\eta_C^t}(y)\leq I(y).\nonumber\end{equation}
Assumption (H3) implies the exponential decay of the memory kernel, hence we have
\begin{equation} I(y)\leq C, \qquad \text{for} \ \ y\geq 1,\nonumber\end{equation}
and as a consequence
\begin{equation}
\lim_{y\rightarrow\infty}I(y)=0,\nonumber
\end{equation}
which yields \eqref{c7}. The lemma is proved.
\end{proof}
 \noindent Lemma \ref{zc} and Lemma \ref{com7} yield the compactness result we
 need
 \bl\label{zc1}
For any $z_0\in \mathcal{V}^0$, $\bigcup_{t\geq 0}z_C(t)$  is
relatively compact in $\mathcal{V}^0$.
 \el
\noindent \textbf{Proof of Theorem \ref{MAIN1}.} On account of Lemma
\ref{abs1}, Lemma \ref{zd}, Lemma \ref{zc1} and the classical result
in dynamical system \cite[Theorem I.1.1]{T}, we can prove the
conclusion of Theorem \ref{MAIN1}, i.e., problem
\eqref{1a}--\eqref{1c} possesses a compact global attractor
$\mathcal {A}$ in $\mathcal{V}^0$.

\section{Convergence to Equilibrium and Convergence Rate} \setcounter{equation}{0}
In this section we prove the convergence of global solutions to
single steady states as time tends to infinity. Let $\mathcal{S}$ be
the set of steady states of $S(t)$,
 \be
 \mathcal{S}=\{Z\in\mathcal{V}^0: S(t)Z=Z,\quad \text{for all}\ \  t\geq
 0\}.
 \ee
It is clear that every steady state $Z_\infty$ has the form
$Z_\infty=(u_\infty,0,0,0)^T$, where $u_\infty$ solves the following
equation
 \be A^2 u_\infty+f(u_\infty)=0,\qquad x\in \Omega, \label{ui}\ee
 with boundary conditions
 \be
 u_\infty=\Delta u_\infty=0,\qquad x\in \Gamma.\label{uii}
 \ee
 The total energy
 \be
 E(t)= \frac{1}{2}\|z(t)\|_{\mathcal{V}^0}^2 + \int_{\Omega} F(u)dx
  \label{Lya}\ee
with $F(u)=\int_0^u f(z)dz$ serves as a Lyapunov functional for
problem \eqref{evo}. Namely, we have
  \be     \frac{d}{dt} E(t) =-\|\nabla
v\|^2+\frac12\int_0^\infty\mu'(s)\|A^{1/2}\eta^t(s)\|^2ds\leq 0,
\qquad \forall\  t>0. \label{Lya1}
 \ee
For any initial data $z_0\in \mathcal{V}^0$, its $\omega$-limit set
is defined as follows:
$$\omega(z_0) = \{z_\infty= (u_\infty,v_\infty,\theta_\infty,\eta_\infty)^T \mid\ \exists \ \{t_n\}\  \text{such that}\
z(t_n) \rightarrow z_\infty\in\ \mathcal{V}^0,\ \text{as}\
t_n\rightarrow +\infty \}.$$ Then we have
 \bl \label{LA}
 \label{cwlas} For any $z_0 \in \mathcal{V}^0$, the
$\omega$-limit set of $z_0$ is a nonempty compact connected subset in $\mathcal{V}^0$. Furthermore, \\
(i) $\omega(z_0)$ is  invariant under the nonlinear semigroup $S(t)$
defined by the solution $z(t,x)$, i.e., $S(t)\omega(z_0) =
\omega(z_0)$ for
all $ t \geq 0$. \\
(ii) $E(t)$ is constant on $\omega(z_0)$. Moreover,
$\omega(z_0)\subset \mathcal{S}$.
 \el
\begin{proof}
Since our system has a continuous Lyapunov functional $E(t)$, the
conclusion of the present lemma follows from Lemma \ref{zc}, Lemma
\ref{zc1} and the well-known results in dynamical system (see, e.g.,
\cite[Lemma~ I.1.1]{T}).
\end{proof}

\br \label{bd3} Since solutions to problem \eqref{ui}\eqref{uii} are
smooth, points in $\omega(z_0)$ are smooth. In particular,
$\omega(z_0)$ is contained in a bounded set in $\mathcal{V}^1$. \er

After the previous preparations, we are ready to finish the proof of
Theorem \ref{MAIN2}.

\noindent \textbf{Part I. Convergence to Equilibrium} \\
For any initial datum $z_0\in \mathcal{V}^0$, it follows from Lemma
\ref{zc} and Lemma \ref{zc1} that there is an equilibrium
$(u_\infty,0,0,0)^T\in \omega(z_0)$ and an increasing unbounded
sequence $\{t_n\}_{n\in\mathbb{N}}$ such that
   \be \lim_{t_n\rightarrow +\infty} (\|u(t_n)-
   u_\infty\|_{V^2}+\|v(t_n)\|+\|\theta(t_n)\|+\|\eta^t(t_n)\|_{\mathcal{M}^1})
   =0. \label{secon}
   \ee
Actually, the convergence for $v,\theta,\eta^t$ can be proved
directly as follows
 \bl Under the assumptions in Theorem \ref{MAIN2}, we
have
\begin{equation}  v(t) \to 0, \quad \theta(t)\to 0,\qquad \text{in} \quad
L^2(\Omega),\label{conv}
\end{equation}
and
\begin{equation}
\eta^t \to 0 \quad \text{in} \quad \mathcal{M}^1,\label{coneta}
\end{equation}
as time goes to infinity.
 \el
 \begin{proof}
Taking the inner product of \eqref{evo} with $z$ in
$\mathcal{V}^{-1}$, we get
 \bea
 \frac12\frac{d}{dt}\|z(t)\|^2_{\mathcal{V}^{-1}}&=&
 -\|v\|^2+\frac12\int_0^\infty\mu'(s)<\eta^t,\eta^t>ds-<f(u),v>_{V^{-1}}\non\\
&\leq&
 -<f(u),v>_{V^{-1}},
 \label{vvv}
 \eea
where in the last step we use (H2). Then the H\"older inequality,
\eqref{w1} and the Sobolev embedding Theorem yield
 \bea
  \frac12\frac{d}{dt}( \|v\|_{V^{-1}}^2+
 \|\theta\|_{V^{-1}}^2+\|\eta^t\|_{\mathcal{M}^{0}}^2)
 &\leq&-<u,v>_{V^{1}}
-<f(u),v>_{V^{-1}}\non\\
 &\leq& C\|Au\|\|v\|+C\|f(u)\|\|v\|\non
 \\
 &\leq &C,\label{vvv2}
 \eea
where $C$ is a constant depending on $\|z_0\|_{\mathcal{V}^0}$.\\
Multiplying \eqref{aa2} by $2\varepsilon$ and adding it to
\eqref{ly} yields
\begin{eqnarray}
&&
\frac{d}{dt}\left(\frac{1}{2}\|z(t)\|^2_{\mathcal{V}^0}+\int_\Omega
F(u)dx+2\varepsilon J(t)\right) +(1-\varepsilon)\|\nabla v\|^2+
\varepsilon\kappa_0\|\theta\|^2\non\\
&&
-\left(\frac12-2C_1\varepsilon\right)\int_0^\infty\mu'(s)\|A^{1/2}\eta^t(s)\|^2ds\leq
2C_2\varepsilon\|\eta^t\|_{\mathcal{M}^1}^2.
\end{eqnarray}
It follows from (H3) that
\begin{eqnarray}
&&
\frac{d}{dt}\left(\frac{1}{2}\|z(t)\|^2_{\mathcal{V}^0}+\int_\Omega
F(u)dx+2\varepsilon J(t)\right) +(1-\varepsilon)\|\nabla v\|^2+
\varepsilon\kappa_0\|\theta\|^2\non\\
&&
+\left[\delta\left(\frac12-2C_1\varepsilon\right)-2C_2\varepsilon\right]
\|\eta^t\|_{\mathcal{M}^1}^2\leq 0.\label{kkkk}
\end{eqnarray}
Taking $\varepsilon$ sufficiently small and integrating \eqref{kkkk}
with respect to $t$, we get
 \be
 \int_0^\infty  \left (\|\nabla v\|^2+\|\theta\|^2+
 \|\eta^t\|_{\mathcal{M}^1}^2\right)dt<\infty.\label{kkk}
 \ee
 Denote
 \be
 h(t)=  \|v\|_{V^{-1}}^2+
 \|\theta\|_{V^{-1}}^2+\|\eta^t\|_{\mathcal{M}^{0}}^2.
  \ee
Then from the continuous embedding $V^1\hookrightarrow
V^0\hookrightarrow V^{-1}$, $\mathcal{M}^1\hookrightarrow
\mathcal{M}^0$, we can conclude from \eqref{kkk} that
 \be h(t)\in L^1(0,\infty).\ee
This and \eqref{vvv2} imply
 \be
 \lim_{t\rightarrow +\infty} h(t)=0. \label{conlo}
  \ee
  Finally, \eqref{conv} and \eqref{coneta} follow from \eqref{secon} and
\eqref{conlo}. The proof is complete.
 \end{proof}

In order to complete the proof of Theorem \ref{MAIN2}, it remains to
show the convergence of $u$. This can be done by making use of a
suitable \L ojasiewicz--Simon type inequality. In our case, it would
be convenient to apply the abstract version in \cite{HJ99}. Denote
 \be
 \mathcal{E}(u)=\frac12\int_\Omega |A u|^2dx + \int_\Omega F(u)dx.
 \ee
 Then we have
 \bl\label{ls}{\rm [\L ojasiewicz--Simon Type Inequality]} Suppose
 that
 assumptions (F1)'(F2) are satisfied.
 Let $\psi$ be a critical point of $\mathcal{E}(u)$. There exist
 constants
 $\rho\in(0,\frac12)$ and $\beta>0$ depending on $\psi$ such that
 for any $u\in V^2$ satisfying $\|u-\psi\|_{V^2}<\beta$,
 there holds
 \be
 \|A^2u+f(u)\|_{V^{-2}}\geq
 |\mathcal{E}(u)-\mathcal{E}(\psi)|^{1-\rho}.
 \ee
 \el
 \br We note that a "smooth" version of \L ojasiewicz--Simon
 inequality of similar type has been introduced in \cite{HR02}. However, the solution to our
problem no longer enjoys the smooth property as in \cite{HR02}.
 \er

We prove the convergence result following a simple argument
introduced in \cite{J981} in which the key observation is that after
certain time $t_0$, the solution $u$ will fall into the small
neighborhood of $u_\infty$ and stay there forever. Unlike parabolic
equations, in order to apply the \L ojasiewicz--Simon approach to
our problem we have to introduce an auxiliary functional which is
usually a perturbation of the Lyapunov functional $E(t)$ due to the
structure of \eqref{evo} (cf. \cite{HJ99,HR02,Mola,WGZ1} and the
references cited therein).

Define
 \bea H(t)&=&\frac{1}{2}\|v(t)\|^2+\frac12\|\theta(t)\|^2+\frac12\|\eta^t\|_{\mathcal{M}^1}^2+
 \mathcal{E}(u(t))-\alpha\int_0^\infty \mu(s)<\theta(t),\eta^t(s)>ds \non\\
 && +\varepsilon <A^2u(t)+f(u(t)),v(t)>_{V^{-2}},\label{Hfun}
 \eea
 where $\alpha>0, \varepsilon>0$ are two coefficients to be
 determined later. It's easy to check that $H(t)$ is well defined for $t\geq 0$. A
direct calculation yields
 \bea
 \frac{dH}{dt}&=& -\|\nabla
v\|^2+\frac12\int_0^\infty\mu'(s)\|A^{1/2}\eta^t(s)\|^2ds-\alpha\int_0^\infty\mu(s)<\Delta
u_t(t),\eta^t(s)>ds \non\\
&&  -\alpha\kappa_0\|\theta\|^2-\alpha\left\|\int_0^\infty\mu(s)
A^\frac12\eta^t(s)ds\right\|^2
+\alpha\int_0^\infty\mu(s)<\theta(t),\eta^t_s(s)>ds\non\\
&& +\varepsilon\left[\|A^2u+f(u)\|^2_{V^{-2}}+ <A^2u+f(u),
Av-A\theta>_{V^{-2}}+ <A^2v+f'(u)v,
v>_{V^{-2}}\right].\non\\\label{dh}
 \eea
It follows from the H\"older inequality, the Poincar\'e inequality
and the Sobolev embedding Theorem that
 \bea \left| <A^2u+f(u),
 Av-A\theta>_{V^{-2}}\right|&\leq& \frac12\|A^2u+f(u)\|^2_{V^{-2}}+
 \|v\|^2+\|\theta\|^2\non\\&
 \leq& \frac12\|A^2u+f(u)\|^2_{V^{-2}}+
 C_9\|\nabla v\|^2+\|\theta\|^2,
 \eea
 \be
 \left|<A^2v+f'(u)v,
 v>_{V^{-2}}\right|\leq \|v\|^2+ \|f'(u)\|_{L^\infty}\|v\|^2\leq
 C_{10} \|\nabla v\|^2,\label{dh1}
 \ee
 where $C_9=C_P^2$ and $C_{10}>0$ depends on $C_P$ and
 $\|z_0\|_{\mathcal{V}^0}$. Recalling \eqref{aa2} and (H3), we deduce from
\eqref{dh}-\eqref{dh1} that
 \bea
 \frac{dH}{dt}&\leq & -\left[1-\frac12\alpha-(C_9+C_{10})\varepsilon\right]\|\nabla
 v\|^2+\left(\frac12-C_1\alpha\right)
 \int_0^\infty\mu'(s)\|A^{1/2}\eta^t(s)\|^2ds\non\\
 && -\left(\frac{\alpha\kappa_0}{2}-\varepsilon\right)\|\theta\|^2+C_2
 \alpha\|\eta^t(s)\|_{\mathcal{M}^1}^2-\frac12\varepsilon\|A^2u+f(u)\|^2_{V^{-2}}\non\\
 &\leq& -\left[1-\frac12\alpha-(C_9+C_{10})\varepsilon\right]\|\nabla
 v\|^2-\left[\left(\frac12-C_1\alpha\right)\delta-C_2\alpha\right]\|\eta^t(s)\|_{\mathcal{M}^1}^2\non\\
 &&
 -\left(\frac{\alpha\kappa_0}{2}-\varepsilon\right)\|\theta\|^2-\frac12\varepsilon\|A^2u+f(u)\|^2_{V^{-2}}.
 \eea
We take $\alpha>0$ small enough such that
 \be
 \left(\frac12-C_1\alpha\right)\delta-C_2\alpha\geq \frac14\delta\quad \text{and}\ \
 \frac12\alpha\leq \frac14,
 \ee
 namely,
 \be
 0<\alpha\leq\min\left\{\frac{\delta}{4(C_1\delta+C_2)}, \frac12\right\}.
  \ee
After fixing $\alpha$, we take $\varepsilon>0$ sufficiently small
satisfying
 \be 0<\varepsilon\leq \min\left\{  \frac{1}{4(C_9+C_{10})},
 \frac14\alpha\kappa_0\right\}.
 \ee
As a result, there exists a positive constant $\gamma$ such that
 \be \frac{d}{dt}H(t)\leq  -\gamma\left(\|\nabla
 v\|^2+\|\eta^t(s)\|_{\mathcal{M}^1}^2+\|\theta\|^2+\|A^2u+f(u)\|^2_{V^{-2}}\right).\label{Lyap}
 \ee
Thus  $H(t)$ is decreasing on $[0,\infty)$. Because $H(t)$ is
bounded from below, it has a finite limit as time goes to infinity.
On the other hand, it follows from \eqref{secon}-\eqref{coneta} that
as $t_n\to \infty$
 \be H(t_n)\to E_\infty=\mathcal{E}(u_\infty).
 \ee From \eqref{Lyap} we can infer that
$H(t) \geq \mathcal{E}(u_\infty)$ for all $t>0$, and the equality
sign holds if and only if $u$ is independent of $t$ and solves
problem (\ref{sta}) while $\theta=v=\eta^t=0$.

We now consider all possibilities.

  \textbf{ Case 1}. If there is a $t_0>0$ such that at this time
   $H(t_0)=\mathcal{E}(u_\infty)$, then for all $t>t_0$, we deduce from
   (\ref{Lyap}) that
   \be \|\nabla v\|\equiv 0. \label{muee}\ee
Namely, $u$ is independent of time for all $t>t_0$. Due to
\eqref{secon}, we can see that (\ref{convu}) holds.

\textbf{Case 2}. For all $t>0$, $H(t)>\mathcal{E}(u_\infty)$. In
this case, there holds
 \be
 -\frac{d}{dt}(H(t)-\mathcal{E}(u_\infty))^{{\rho}}=-{\rho}
   (H(t)-\mathcal{E}(u_\infty))^{{\rho}-1}\frac{d}{dt}H(t), \label{5dd}
 \ee
 here $\rho\in (0,\frac12)$ is the exponent in Lemma \ref{ls}. By the H\"older inequality, we obtain
 \bea (H -\mathcal{E}(u_\infty))^{1-{\rho}}&\leq& C\left(\|v\|^{2(1-{\rho})}+\|\theta\|^{2(1-{\rho})}
 +\|\eta^t\|_{\mathcal{M}^1}^{2(1-{\rho})}+|\mathcal{E}(u)-\mathcal{E}(u_\infty)|^{1-{\rho}}\right.\non\\
 && \left.\ \ +
 \|\theta\|^{1-{\rho}}\|\eta^t\|_{\mathcal{M}^1}^{1-{\rho}}+\|A^2u+f(u)\|^{1-{\rho}}_{V^{-2}}\|v\|^{1-{\rho}}\right).
 \eea
 Besides, the Young inequality yields
 \be
 \|A^2u+f(u)\|^{1-{\rho}}_{V^{-2}}\|v\|^{1-{\rho}}\leq
 \|A^2u+f(u)\|_{V^{-2}}+ C \|v\|^{\frac{1-{\rho}}{\rho}}.
 \ee
Noting that $(1-{\rho})/{\rho}>1$ and $2(1-{\rho})>1$, by the
uniform bounds obtained in previous section, we conclude
 \bea (H
 -\mathcal{E}(u_\infty))^{1-{\rho}}&\leq& C\left(\|v\|+\|\theta\|
 +\|\eta^t\|_{\mathcal{M}^1}+\|A^2u+f(u)\|_{V^{-2}}+|\mathcal{E}(u)-\mathcal{E}(u_\infty)|^{1-{\rho}}\right).\non\\
 \label{dh2}
 \eea
 It follows from \eqref{secon} that there exists $N\in \mathbb{N}$
 such that for any $n\geq N$, $\|u(t_n)-u_\infty\|_{V^2}<\beta$. Set
 \be \bar{t}_n = \sup \left\{\ t>t_n \left\vert \ \|u(\tau)-u_\infty\|_{V^2}< \beta,\   \forall \tau\in [t_n , t] \right.\right\}.
 \label{4.39dual} \ee
 Observe that $\bar{t}_n > t_n$ for all $n \geq
 N$, due to the continuity of the orbit in
$\mathcal{V}^0$. Now we have to deal with two subcases.

(a) There exists $n_0$ such that $\bar{t}_{n_0}=\infty$. By Lemma
\ref{ls}, \eqref{dh}\eqref{5dd}\eqref{dh2} and the Poincar\'e
inequality, we can conclude that
  \bea &&
  -\frac{d}{dt}(H(t)-\mathcal{E}(u_\infty))^{{\rho}}\non\\
  &\geq &
  C\rho\gamma \frac {\|\nabla
 v\|^2+\|\eta^t(s)\|_{\mathcal{M}^1}^2+\|\theta\|^2+\|A^2u+f(u)\|^2_{V^{-2}}}{\|v\|+\|\theta\|
 +\|\eta^t\|_{\mathcal{M}^1}+\|A^2u+f(u)\|_{V^{-2}}+|\mathcal{E}(u)-\mathcal{E}(u_\infty)|^{1-{\rho}}}\non\\
 &\geq& C_{11} \left(\|\nabla
 v\|+\|\eta^t(s)\|_{\mathcal{M}^1}+\|\theta\|+\|A^2u+f(u)\|_{V^{-2}}\right).\label{dh3}
  \eea
Integrating from $t_{n_0}$ to $t$, we obtain
 \bea &&  (H(t)-\mathcal{E}(u_\infty))^{{\rho}}\non\\
 &&\ \ \ + C_{11} \int_{t_{n_0}}^t \left(\|\nabla
 v\|+\|\eta^t(s)\|_{\mathcal{M}^1}+\|\theta\|+\|A^2u+f(u)\|_{V^{-2}}\right)d\tau\non\\
 &\leq& (H(t_{n_0})-\mathcal{E}(u_\infty))^{{\rho}}.
 \eea
Recalling that $H(t) - \mathcal{E}(u_\infty) \geq 0$ for $t>0$, we
infer \be \int_{t_{n_0}}^t \|
     v(\tau)\|_{V^1}  d\tau < \infty, \qquad \forall\  t\geq t_{n_0}. \label{conu3}\ee
Thus, $u(t)$ converges in $V^1$. Then by the precompactness property
of $u(t)$ in $V^2$ (see Section 4), we can conclude \eqref{convu}.

(b) For all $n\in\mathbb{N}$, $\bar{t}_n < \infty$.

Since $H(t)$ is decreasing in $[0, \infty)$ and it has a finite
limit $E_\infty=\mathcal{E}(u_\infty)$ as $t\to \infty$, then for
any $\zeta\in (0,\beta)$ there exists an integer $N$ such that when
$n \geq N$, for all $t \geq t_n>0$, there holds
 \be
     (H(t_n) - \mathcal{E}(u_\infty))^{{\rho}}
      -(H(t) - \mathcal{E}(u_\infty))^{{\rho}}< \frac{C_{11}}{2}\zeta.
     \label{4.38dual}
 \ee
 As a result, for $n\geq N$ there holds
 \be
    \int_{t_n}^{\bar{t}_n} \|v(\tau)\|_{V^1} d\tau < \frac{\zeta}{2}.
 \ee
Moreover, by choosing $N$ sufficiently large we have
 \be
   \|u(t_n)-u_\infty\|_{V^2}< \frac{\zeta}{2}, \quad\forall\, n\geq N.
\label{4.37dual}
 \ee
 These imply that
 \be \|u(\bar{t}_n) - u_\infty \|_{V^1}
 \leq  \| u(t_n) - u_\infty\|_{V^1} + \int_{t_n}^{\bar{t}_n}
 \|v(\tau)\|_{V^1} d\tau < \zeta,\quad \forall\ n\geq N.
 \ee
  Therefore,
 \be \lim_{\bar{t}_n\rightarrow +\infty} \| u(\bar{t}_n) -
 u_\infty\|_{V^1}=0.
 \ee
 On the other hand, the precompactness of $u$ in $V^2$ implies that there exists a
subsequence of $\{u(\bar{t}_n)\}$, still denoted by
$\{u(\bar{t}_n)\}$, converging to $u_\infty$ in $V^2$. Thus for $n$
sufficiently large, we get \be \|u(\bar{t}_n) - u_\infty \|_{V^2} <
\beta \ee which contradicts the definition of $\bar{t}_n$.\bigskip

 \textbf{Part II.
Convergence Rate}\\
For $t\geq t_0$ with $t_0$ sufficiently large, it follows from Lemma
\ref{ls} and \eqref{dh2}(\ref{dh3}) that
 \be \frac{d}{dt}(H(t) - \mathcal{E}(u_\infty))+ C\left(H(t) -
 \mathcal{E}(u_\infty)\right)^{2(1-\rho)}\ \leq\ 0. \label{rate1}
 \ee
This yields (cf. \cite{HJ99,Z,WGZ1})
 \be H(t)-\mathcal{E}(u_\infty) \leq  C (1+t)^{-1/(1-2\rho)},
 \qquad\forall\,t\geq t_0.
 \ee
 Integrating \eqref{dh2} on
$(t,\infty)$, we have
 \be \int_{t}^{\infty}
 \|v\|_{V^1} d\tau \leq C
 (1+t)^{-\rho/(1-2\rho)},\qquad\forall\,t\geq t_0.\label{rate2}
 \ee
 By adjusting the constant $C$ properly, we obtain
 \be
    \|u(t)-u_\infty\|_{V^1}\leq C(1+t)^{-\rho/(1-2\rho)}, \quad\forall\ t\geq 0.\label{rate3}
 \ee
 Based on this estimate for $u$ in $V^1$ norm, we are able to obtain
 the estimates (in higher order norm) stated in Theorem \ref{MAIN2}.

By subtracting the evolution equations \eqref{1a1} and their
corresponding stationary equations \eqref{sta}, we have
\begin{equation}
  \left\{\begin{array}{l} \theta_t -\Delta u_t -\int_0^\infty \mu(s)\Delta \eta^t(s)ds =0, \\
    u_{tt} -\Delta u_t+ \Delta\theta +\Delta^2 (u-u_\infty) +f(u)- f(u_\infty) =0,
   \end{array}
   \right.\label{1aa}
 \end{equation}
Similar to \eqref{ly}, we can see that
\begin{eqnarray}
&&\frac{d}{dt}\left(\frac{1}{2}\|u(t)-u_\infty\|^2_{V^2}
+\frac12\|v\|^2+\frac12 \|\theta\|^2+\frac12
\|\eta^t\|_{\mathcal{M}^1}^2+\int_\Omega F(u)dx\right.\non\\
&& \left.-\int_\Omega
F(u_\infty)dx- \int_\Omega f(u_\infty)(u-u_\infty)dx  \right)\non\\
&&  +\|\nabla
v\|^2-\frac12\int_0^\infty\mu'(s)\|A^{1/2}\eta^t(s)\|^2ds\non\\
&=& 0.\label{lya}
\end{eqnarray}
Multiplying the second equation in \eqref{1aa} by $u-u_\infty$ and
integrating on $\Omega$, we get
\begin{eqnarray}
&& \frac{d}{dt}\left( \frac12\|\nabla u-\nabla u_\infty\|^2+
\int_\Omega v(u-u_\infty)dx\right)- \|v\|^2+ \|A(u-u_\infty)\|^2\non
\\
&=& -\int_\Omega (f(u)-f(u_\infty))(u-u_\infty)dx +\int_\Omega
\theta A(u-u_\infty)dx. \label{dbb}
\end{eqnarray}
Multiplying \eqref{aa2} by $2\varepsilon$ and multiplying
\eqref{dbb} by $\varepsilon^2$ respectively, then adding the
resultants to \eqref{lya} yield
\begin{eqnarray}
&& \frac{d}{dt}\left( \frac{1}{2}\|u(t)-u_\infty\|^2_{V^2}
+\frac12\|v\|^2+\frac12 \|\theta\|^2+\frac12
\|\eta^t\|_{\mathcal{M}^1}^2+\frac{\varepsilon^2}{2}\|\nabla
(u-u_\infty)\|^2+\int_\Omega F(u)dx\right.\non\\
&& \left.-\int_\Omega F(u_\infty)dx - \int_\Omega
f(u_\infty)(u-u_\infty)dx
+ \varepsilon^2\int_\Omega v(u-u_\infty) dx+2\varepsilon J(t)\right)\non\\
& & \ +(1-\varepsilon)\|\nabla v\|^2 +\varepsilon\kappa_0
\|\theta\|^2+\left[\left(\frac12-2C_1\varepsilon\right)\delta-2C_2\varepsilon\right]\|\eta^t\|_{\mathcal{M}^1}^2+
\varepsilon^2\|A(u-u_\infty)\|^2\nonumber\\
&\leq& \varepsilon^2\|v\|^2 -\varepsilon^2\int_\Omega
(f(u)-f(u_\infty))(u-u_\infty)dx+\varepsilon^2\int_\Omega \theta A
(u-u_\infty)dx. \label{rate4}
\end{eqnarray}
We now estimate the three terms on the right hand side of inequality
\eqref{rate4}.
 \be \varepsilon^2\|v\|^2\leq \varepsilon^2C_P^2\|\nabla v\|^2,
 \ee
 \bea \left|-\varepsilon^2\int_\Omega
 (f(u)-f(u_\infty))(u-u_\infty)dx\right|\leq
 \varepsilon^2\|f'\|_{L^\infty}\|u-u_\infty\|^2
         \leq  C\varepsilon^2\|u-u_\infty\|^2, \label{scr4}
 \eea
 \bea
 \left|\varepsilon^2\int_\Omega \theta A
(u-u_\infty)dx\right|\leq \frac14\varepsilon^2 \|A (u-u_\infty)\|^2+
\varepsilon^2\|\theta\|^2.\label{rate4aa}
 \eea
On the other hand, by the Taylor's expansion, we have \be
F(u)=F(u_\infty) + f(u_\infty)(u-u_\infty) + f'(\xi)(u-u_\infty)^2,
\ee
where $\xi=a u+ (1-a)u_\infty$ with $a\in [0,1]$.\\
 Then we deduce that
\bea & & \left\vert\int_\Omega F(u)dx-\int_\Omega
        F(u_\infty)dx+\int_\Omega f(u_\infty)u_\infty
         dx-\int_\Omega f(u_\infty)u dx\right\vert\non\\
&  =  & \left\vert \int_\Omega f'(\xi)(u-u_\infty)^2 dx\right\vert\non\\
& \leq &
          \|f'(\xi)\|_{L^\infty}\|u-u_\infty\|^2\leq C\|u-u_\infty\|^2.
          \label{scr10}
  \eea
  Let us define now, for $t\geq 0$,
 \bea
 y(t)&=& \frac{1}{2}\|u(t)-u_\infty\|^2_{V^2}
+\frac12\|v\|^2+\frac12 \|\theta\|^2+\frac12
\|\eta^t\|_{\mathcal{M}^1}^2+\frac{\varepsilon^2}{2}\|\nabla
(u-u_\infty)\|^2+\int_\Omega F(u)dx\non\\
&& -\int_\Omega F(u_\infty)dx - \int_\Omega
f(u_\infty)(u-u_\infty)dx + \varepsilon^2\int_\Omega v(u-u_\infty)
dx+2\varepsilon J(t).
 \eea
Taking $\varepsilon$ sufficiently small, it follows from the
H\"older inequality, \eqref{scr10} and \eqref{JJ} that there exist
constants $\gamma_0,\gamma_1,\gamma_2>0$ such that
 \be
 \gamma_0\|z-z_\infty\|^2_{\mathcal{V}^0}\geq y(t)\geq
 \gamma_2\|z-z_\infty\|^2_{\mathcal{V}^0}- \gamma_1\|u-u_\infty\|_{V^1}^2.\label{yy}
 \ee
Moreover, for small $\varepsilon$ we can deduce from
\eqref{rate4}--\eqref{rate4aa} and \eqref{yy} that for certain
$\gamma_3>0$, the following inequality holds
 \be \frac{d}{dt}y(t)+ \gamma_3 y(t)\leq C\|u-u_\infty\|_{V^1}^2.\ee
 The Gronwall inequality and \eqref{rate3} yield (see e.g.,
 \cite{W07,WGZ1})
 \be y(t)\leq C(1+t)^{-2\rho/(1-2\rho)}, \quad \forall \ t\geq 0,\ee
 which together with \eqref{yy} implies that
 \be
 \|z-z_\infty\|_{\mathcal{V}^0}\leq  C(1+t)^{-\rho/(1-2\rho)}, \quad \forall \ t\geq
 0.\label{zrate}
 \ee
 The proof of Theorem \ref{MAIN2} is now complete.

Before ending this paper, we give a further remark on the estimate
of convergence rate. As has been shown in the previous section, the
solution $z(t)$ to our problem \eqref{evo} with initial data $z_0\in
\mathcal{V}^0$ can be decomposed into two parts
$z(t)=z_D(t)+z_C(t),$ where
$z_D(t)=\left(u_D(t),v_D(t),\theta_D(t),\eta_D^t\right)^T$ and
$z_C(t)=\left(u_C(t),v_C(t) ,\theta_C(t),\eta_C^t\right)^T$ satisfy
\eqref{evoD} and \eqref{evoC} respectively. It is also shown in
Lemma \ref{zd} that $z_D(t)$ will  decay to $0$ in $\mathcal{V}^0$
exponentially fast. This convergence rate is obviously better than
the rate for  $z(t)$ obtained in Theorem \ref{MAIN2}. As a result,
we can easily obtain the following result for the compact part
$z_C(t)$ from Lemma \ref{zd} and Theorem \ref{MAIN2}:

\begin{proposition}\label{zc1a}
Under the assumptions of Theorem \ref{MAIN2}, we have
 \be \|z_C(t)-z_\infty\|_{\mathcal{V}^0}\leq C(1+t)^{-\rho/(1-2\rho)}, \quad \forall \ t\geq
 0,\label{ratezc1}\ee
where $C>0$ is a constant depending on $\|z_0\|_{\mathcal{V}^0}$ and
$z_\infty=(u_\infty,0,0,0)^T$.
\end{proposition}
\begin{proof}
We notice that
 \be \|z_C(t)-z_\infty\|_{\mathcal{V}^0}\leq
 \|z(t)-z_\infty\|_{\mathcal{V}^0}+\|z_D(t)\|_{\mathcal{V}^0},\ee
 \be \lim_{t\to +\infty} e^{-\frac{\delta_1}{2}t}(1+t)^{\rho/(1-2\rho)}=0.\ee
Then the  conclusion \eqref{ratezc1} follows from Lemma \ref{zd} and
\eqref{zrate} after the constant $C$ is properly modified.
\end{proof}
Moreover, Lemma \ref{zc} provides a uniform estimate of $z_C$ in
$\mathcal{V}^1$. As a direct consequence, this fact and Proposition
\ref{zc1a} imply the weak convergence of $z_C$ such that
$$ z_C(t)\rightharpoonup z_\infty,\quad \text{in}\ \mathcal{V}^1, \quad
\text{as}\ t\to +\infty.$$
 Based on the idea we used in the proof of Theorem \ref{MAIN2},
 we are able to get a stronger result, namely
 \bt
Under the assumptions of Theorem \ref{MAIN2}, we have
 \be \|z_C(t)-z_\infty\|_{\mathcal{V}^1}\leq C(1+t)^{-\rho/(1-2\rho)}, \quad \forall\ t\geq
 0,\label{ratezc}\ee
where $C>0$ is a constant depending on $\|z_0\|_{\mathcal{V}^0}$ and
$\|u_\infty\|_{V^3}$.
 \et
\begin{proof}
Subtracting \eqref{sta} from \eqref{evoC} we have
\begin{equation}
  \left\{\begin{array}{l} \displaystyle \frac{d}{dt}(z_C-z_\infty)=L(z_C-z_\infty)+(0,-f(u)+f(u_\infty), 0, 0 )^T, \\
    (z_C-z_\infty)|_{t=0}=z_\infty.
   \end{array}
   \right.\label{evoC1}
 \end{equation}
Taking the inner product of the resulting system \eqref{evoC1} and
$A(z_C-z_\infty)$ in $ \mathcal{V}^0$, we get
 \bea
&&
\frac{d}{dt}\left(\frac12\|A^{1/2}(z_C-z_\infty)\|^2_{\mathcal{V}^0}+\int_\Omega
(f(u)-f(u_\infty))A(u_C-u_\infty)dx\right)\non\\
&& +\|A
v_C\|^2 -\frac12\int_0^\infty\mu'(s)\|A\eta_C^t(s)\|^2ds\non\\
&=&\int_\Omega f'(u)vA(u_C-u_\infty).\label{COM1a}
 \eea
Next, multiplying the second equation in \eqref{evoC1} by
$A(u_C-u_\infty)$ and integrating on $\Omega$, we have
 \bea
&& \frac{d}{dt}\left(\int_\Omega v_CA(u_C-u_\infty)
dx+\frac{1}{2}\|A(u_C-u_\infty)\|^2
\right)+\|A^{3/2}(u_C-u_\infty)\|^2-\|A^{1/2} v_C\|^2 \non\\
&=& - \int_\Omega (f(u)-f(u_\infty))A(u_C-u_\infty) dx+\int_\Omega
A\theta_C A(u_C-u_\infty)dx.\label{COM2a}
 \eea
 Now we introduce the functional
 \bea \Upsilon(t)&=& \|A^{1/2}(z_C(t)-z_\infty)\|^2_{\mathcal{V}^0}+
 2\int_\Omega (f(u)-f(u_\infty))A(u_C-u_\infty) dx\non\\
 &&+ 2\varepsilon^2\int_\Omega v_CA(u_C-u_\infty)dx  +
 \varepsilon^2\|A (u_C-u_\infty)\|^2+4\varepsilon J_C(t).
 \eea
It follows from Theorem \ref{MAIN2}, Proposition \ref{zc1a} and
\eqref{JC} that for $t\geq 0$,
 \bea
 && \left|\int_\Omega (f(u)-f(u_\infty))A(u_C-u_\infty) dx\right|\non\\
 &\leq&
 \| f'\|_{L^\infty}\|u-u_\infty\|\|A(u_C-u_\infty)\|\leq
 C(1+t)^{-2\rho/(1-2\rho)},
 \eea
 \be
 \left| \int_\Omega v_CA(u_C-u_\infty)dx\right|\leq
 \|v_C\|\|A(u_C-u_\infty)\|\leq C(1+t)^{-2\rho/(1-2\rho)},
 \ee
 \be
 |J_C(t)|\leq C\|z_C-z_\infty\|^2_{\mathcal{V}^1}.
 \ee
 As a result, after choosing $ \varepsilon >0$ sufficiently small, there is a constant $C>0$ such that
\begin{equation}
\|z_C(t)-z_\infty\|^2_{\mathcal{V}^1}\leq 2\Upsilon(t)+
C(1+t)^{-2\rho/(1-2\rho)}.\label{upsilon}
\end{equation}
It follows from \eqref{COM1a}--\eqref{COM2a} and \eqref{COM3} that
 \bea
 &&\frac12\frac{d}{dt}\Upsilon(t)+ (1-\varepsilon)\|A
v_C\|^2-\left(\frac12-2C_7\varepsilon\right)\int_0^\infty\mu'(s)\|A\eta_C^t(s)\|^2ds
\non\\
&&-\varepsilon^2\|A^{1/2}v_C\|^2+\varepsilon^2\|A^{3/2}(u_C-u_\infty)\|^2
+\varepsilon
\kappa_0\|A^{1/2}\theta_C\|^2\non\\
&\leq& \int_\Omega f'(u)vA(u_C-u_\infty)dx
+2C_6\varepsilon\|\eta_C^t\|_{\mathcal{M}^2}^2\non\\
&& -\varepsilon^2 \int_\Omega
(f(u)-f(u_\infty))A(u_C-u_\infty)dx+\varepsilon^2\int_\Omega
A\theta_CA(u_C-u_\infty)dx.\label{c1a}
 \eea
 The right hand side of \eqref{c1a} can be estimated as follows
\begin{eqnarray}
\left| \int_\Omega f'(u)v A(u_C-u_\infty)dx\right| \leq
\|f'(u)\|_{L^\infty}\|v\|\|A(u_C-u_\infty)\|\leq
C(1+t)^{-2\rho/(1-2\rho)}, \label{ff5}
\end{eqnarray}
 \bea&& \left|\int_\Omega
(f(u)-f(u_\infty))A(u_C-u_\infty)dx\right|\non\\
&\leq& \|f'\|_{L^\infty}\|u-u_\infty\|\|A(u_C-u_\infty)\|\leq
 C(1+t)^{-2\rho/(1-2\rho)},
 \eea
\bea
 \left| \varepsilon^2\int_\Omega
A\theta_CA(u_C-u_\infty)dx\right|\leq
\frac{\varepsilon\kappa_0}{2}\|A^{1/2}\theta_C\|^2+C\varepsilon^3\|A^{3/2}(u_C-u_\infty)\|^2,\label{c2a}
\eea where in the above estimates $C$ is a constant depending on
$\|z_0\|_{\mathcal{V}^0}$ at most. Similar to the previous section,
we can choose $\varepsilon>0$ small enough and consequently there is
a constant $\gamma_4>0$ such that
 \be
 \frac{d}{dt}\Upsilon(t)+ \gamma_4\Upsilon(t)\leq C(1+t)^{-2\rho/(1-2\rho)}.
 \ee
As a result,
 \be \Upsilon(t)\leq C(1+t)^{-2\rho/(1-2\rho)}, \quad \forall\ t\geq
 0,\label{up1}\ee
 here $C>0$ is a constant depending on $\|z_0\|_{\mathcal{V}^0}$
 and $\|u_\infty\|_{V^3}$ (see Remark \ref{bd3}).

 The required estimate \eqref{ratezc} follows from \eqref{up1} and
 \eqref{upsilon}. We complete the proof.
\end{proof}

%%%%%%%%%%%%%%%%%%%%%%%%%%%%%%%%%%%%%%%%%%%%%%%%%%%%%%%%%%%%%%%%%%%
\bigskip
\noindent\textbf{Acknowledgments.} The author wants to thank Prof.
S. Zheng for his enthusiastic help and encouragement. The author
also thanks Dr. G. Mola for sharing the interesting results in his
PhD thesis. The research of the author was supported by China
Postdoctoral Science Foundation.
%%%%%%%%%%%%%%%%%%%%%%%%%%%%% References %%%%%%%%%%%%%%%%%%%%%%%%%%


\begin{thebibliography}{99}
\itemsep=0pt

\bibitem{AF} S. Aizicovici and E. Feireisl, Long-time stabilization of solutions
to a phase--field model with memory, J. Evol. Equ., \textbf{1
}(2001), 69--84.

\bibitem{AP} S. Aizicovici and H. Petzeltov\'a, Convergence of solutions of
phase--field systems with a nonconstant latent heat,  Dynam. Systems
Appl.,  \textbf{14}  (2005), 163--173.

\bibitem{BP} V. Belleri and V. Pata, Attractors for semilinear strongly
damped wave equation on $\mathbb{R}^3$, Disc. Contin. Dynam. Sys. A,
\textbf{7} (2001), 719--735.

\bibitem{CMP} V.V. Chepyzhov, E. Mainini and V. Pata,
Stability of abstract linear semigroups arising from heat conduction
with memory, Asymptot. Anal., \textbf{50} (2006), 269--291.


\bibitem{CEL}I. Chueshov, M. Eller and I. Lasiecka, On the attractor
for a semilinear wave equation with critical exponent and nonlinear
boundary dissipation, Comm. PDEs., {\bf 27} (2002), 1901--1951.

\bibitem{CG} B.D. Coleman and M.E. Gutin, Equipresence and constitutive equations
for rigid heat conductors, Z. Angew. Math. Phys., \textbf{18}
(1967), 199--208.

\bibitem{D} C.M. Dafermos, Asymptotic stability in viscoelasticity, Arch. Ration. Mech.
Anal., \textbf{37} (1970), 297--308.

\bibitem{FLR} M. Fabrizio, B. Lazzari and J.E. Mu\~noz Rivera, Asymptotic behavior of
thermoelastic plates of weakly hyperbolic type, Differential
Integral Equations, \textbf{13} (2000), 1347--1370.

\bibitem{GW08} C.G. Gal and H. Wu, Asymptotic behavior of a
Cahn-Hilliard equation with Wentzell boundary conditions and mass
conversation, Disc. Contin. Dynam. Sys. A, to appear.

\bibitem{GGP} S. Gatti, M. Grasselli and V. Pata, Exponential attractors for a conserved phase--field system
with memory, Phys. D, \textbf{189} (2004), 31--48.

\bibitem{GGP1} C. Giorgi, M. Grasselli and V. Pata, Uniform attractors for a phase--field
model with memory and quadratic nonlinearity,  Indiana Univ. Math.
J.,  \textbf{48} (1999), 1395--1445.

\bibitem{GP} C. Giorgi and V. Pata, Stability of linear thermoelastic systems with memory,
Math. Mod. Meth. Appl. Sci., \textbf{11} (2001), 627--644.

\bibitem{GAA} M. Grasselli, On the large time behavior of a phase--field
system with memory, Asymptot. Anal., \textbf{56} (2008), 229--249.

\bibitem{GRP} M. Grasselli, J.E. Mu\~noz Rivera and V. Pata, On the energy decay of
the linear thermoelastic plate with memory, J. Math. Anal. Appl.,
\textbf{309} (2005), 1--14.

\bibitem{GHS} M. Grasselli, H. Petzeltov\'a and G. Schimperna, Asymptotic behavior of a nonisothermal viscous
Cahn--Hilliard equation with inertial term, J. Differential
Equations, \textbf{239} (2007), 38--60.

\bibitem{GP01} M. Grasselli and V. Pata, Upper semicontinuous attractor for a hyperbolic
phase--field model with memory, Indiana Univ. Math. J., \textbf{50}
(2001), 1281--1308.

\bibitem{GPV} M. Grasselli, V. Pata and F.M. Vegni, Longterm dynamics of a
conserved phase--field model with memory, Asymptot. Anal.,
\textbf{33} (2003), 261--320.

\bibitem{GP2}
M. Grasselli and V. Pata, Asymptotic behavior of a
parabolic--hyperbolic system, Commun. Pure Appl. Anal., \textbf{3}
(2004), 849--881.

\bibitem{Gra} M. Grasselli and V. Pata, Existence of a universal attractor for a
fully hyperbolic phase--field system, J. Evol. Equ., \textbf{4}
(2004), 27--51.

\bibitem{WGZ3} M. Grasselli, H. Wu and S. Zheng, Asymptotic behavior of
a non-isothermal Ginzburg--Landau model, Quart. Appl. Math., to
appear.

\bibitem {GP1} M.E. Gurtin and A.C. Pipkin, A general theory of heat conduction with
finite wave speeds, Arch. Ration. Mech. Anal., \textbf{31} (1968),
113--126.

\bibitem{HJ99} A. Haraux and M.A. Jendoubi, Convergence of bounded
weak solutions of the wave equation with dissipation and analytic
nonlinearity, Calc. Var., {\bf 9} (1999), 95--124.

\bibitem{HJ01} A. Haraux and M.A. Jendoubi, Decay estimates to equilibrium for
some evolution equations with an analytic nonlinearity, Asymptot.
Anal., \textbf{26} (2001), 21--36.

\bibitem{He} L. Herrera and D. P\'avon, Hyperbolic theories of dissipation: Why
and when do we need them? Phys. A, \textbf{307} (2002), 121--130.

\bibitem{HR02} K.-H. Hoffmann and P. Rybka, Analyticity of the nonlinear term forces
convergence of solutions to two equations of continuum mechanics,
Nonlinear Analysis, \textbf{50} (2002), 409--424.

\bibitem{J981}M.A. Jendoubi, A simple unified approach to
some convergence theorem of L. Simon, J. Func. Anal., {\bf 153}
(1998), 187--202.

\bibitem{JJ} J. Jiang, Convergence to equilibrium for a parabolic--hyperbolic phase--field
model with Cattaneo heat flux law, J. Math. Anal. Appl.,
\textbf{341} (2008), 149--169.

\bibitem{JJ1} J. Jiang, Convergence to equilibrium for a fully hyperbolic
phase--field model with Cattaneo heat flux law, preprint, 2008.

\bibitem{Jo1} D.D. Joseph and L. Preziosi, Heat waves, Rev. Modern Phys., \textbf{61}
(1989), 41--73.

\bibitem{Jo2} D.D. Joseph and L. Preziosi, Addendum to the paper: Heat waves,
Rev. Modern Phys., \textbf{62} (1990), 375--391.

\bibitem{L} J. Lagnese, Boundary stabilization of thin plates, SIAM Stud. Appl.
 Math. \textbf{10}, SIAM, Philadelphia, 1989.

\bibitem{LZ} Z. Liu and S. Zheng, Semigroups associated with dissipative systems,
Chapman \& Hall/CRC Press, Boca Raton, FL, 1999.

\bibitem{Mola} G. Mola, Convergence to equilibria for a three--dimensional conserved
phase--field system with memory, Elec. J. Diff. Eqns., \textbf{23}
(2008), 1--16.

\bibitem{Molaphd} G. Mola, \textit{Global and Exponential Attractors for a
Conserved Phase--Field System with Grutin--Pinkin Heat Conduction
Law}, Ph.D. thesis, Dipartimento di Matematica ``F. Brioschi'',
Politecnico di Milano, 2007.



\bibitem{RR1} J.E. Mu\~noz Rivera and R. Racke, Smoothing properties, decay and
global existence of solutions to nonlinear coupled systems of
thermoelastic type, SIAM J. Math. Anal., \textbf{26} (1995),
1547--1563.

\bibitem{RR} J.E. Mu\~noz Rivera and R. Racke, Large solutions and smoothing
properties for nonlinear thermoelastic systems, J. Differential
Equations., \textbf{127} (1996), 454--483.

\bibitem{pa} V. Pata, Exponential stability in linear viscoelasticity, Quart. Appl. Math., \textbf{64} (2006), 499--513.

\bibitem{T}{\au R. Temam}, {\ti Infinite-dimensional Dynamical Systems in
Mechanics and Physics,} {\jou Appl. Math. Sci.,} \textbf{68}, Springer-Verlag, New York, 1988.

\bibitem{W07} H. Wu, Convergence to equilibrium for a Cahn--Hilliard
model with the Wentzell boundary condition, Asymptot. Anal.,
\textbf{54} (2007), 71--92.

\bibitem{WGZ1} H. Wu, M. Grasselli and S. Zheng, Convergence to
equilibrium for a parabolic--hyperbolic phase--field system with
Neumann boundary conditions, Math. Mod. Meth. Appl. Sci.,
\textbf{17} (2007), 1--29.

\bibitem{WGZ2} H. Wu, M. Grasselli and S. Zheng, Convergence to
equilibrium for a nonlinear parabolic--hyperbolic phase--field
system with dynamic boundary condition, J.
 Math. Anal. Appl., \textbf{329} (2007), 948--976.


\bibitem{Z}
{\au S.~Zheng}, {\bk Nonlinear Evolution Equations}, Pitman series
Monographs and Survey in Pure and Applied Mathematics, \textbf{133},
\eds{Chapman \& Hall/CRC}{Boca Raton, Florida}{2004}

\end{thebibliography}
\end{document}